\documentclass[11pt]{article}
\usepackage{amsmath,epsfig,subfigure,amsthm,amsfonts,amsbsy,amssymb,latexsym,amsxtra}
\usepackage[usenames]{color}
\usepackage{pst-infixplot}
\usepackage{pst-plot}
\usepackage{pstricks}

\usepackage[sort,comma]{natbib}
\makeatletter \@addtoreset{equation}{section}

\newcommand{\beq}{\begin{equation}}
\newcommand{\eeq}{\end{equation}}
\newcommand{\bed}{\begin{displaymath}}
\newcommand{\eed}{\end{displaymath}}
\newcommand{\ben}{\begin{eqnarray*}}
\newcommand{\een}{\end{eqnarray*}}
\newcommand{\bedd}{\bed\begin{array}{l}}
\newcommand{\eedd}{\end{array}\eed}
\newcommand{\bea}{\bed\begin{array}{rl}}
\newcommand{\eea}{\end{array}\eed}
\newcommand{\ad}{&\!\!\!\disp}
\newcommand{\aad}{&\disp}
\newcommand{\barray}{\begin{array}{ll}}
\newcommand{\earray}{\end{array}}
\newcommand{\be}{\begin{eqnarray}}
\newcommand{\ee}{\end{eqnarray}}
\newcommand{\by}{\begin{eqnarray*}}
\newcommand{\ey}{\end{eqnarray*}}

\def\F{{\mathcal F}}
\def\op{{\mathcal A}}
\def\cd{(\cdot)}

\def\rr{{\mathbb R}}

\newcommand{\la}{\lambda}

\newcommand{\dl}{\delta}

\newcommand{\vphi}{{\varphi}}

\newcommand{\al}{\alpha}

\newcommand{\ex}{{\mathrm E}}

\newcommand{\wrt}{{with respect to }}

\def\({\left(}
\def\){\right)}

\newcommand{\disp}{\displaystyle}

\def\one{{\hbox{1{\kern -0.35em}1}}}
\newcommand{\set}[1]{\left\{#1\right\}}

\newcommand{\df}{\, \mathrm{d}}
\newcommand{\td}{\tilde}

\newtheorem{thm}{Theorem}[section]
\newtheorem{prop}[thm]{Proposition}
\newtheorem{lem}[thm]{Lemma}

\theoremstyle{definition}
\newtheorem{rem}[thm]{Remark}

\newcommand{\thmref}[1]{Theorem~{\rm \ref{#1}}}
\newcommand{\lemref}[1]{Lemma~{\rm \ref{#1}}}

\title{Optimal Dividend Policies for Piecewise-Deterministic Compound Poisson Risk Models\thanks{This research was supported in part by National Science Foundation
 under grant DMS-1108782 and from the UWM Research Growth Initiative, and City University
 of Hong Kong (SRG) 7002677.
}}

\author{ Runhuan Feng, University of Illinois at Urbana-Champaign\\
Hans W. Volkmer, University of Wisconsin - Milwaukee \\
 Shuaiqi Zhang, Hebei University of Technology\\
Chao Zhu, University of Wisconsin - Milwaukee 
}
\begin{document}

\maketitle

\vspace{-1cm}
\begin{abstract}
This paper considers the optimal dividend payment problem in piecewise-deterministic compound Poisson risk models. The objective is to maximize the  expected discounted dividend payout up to the time of ruin. We provide a comparative study in this general framework of both restricted and unrestricted payment schemes, which were only previously treated separately in certain special cases of risk models in the literature. In the case of restricted payment scheme, the value function is shown to be a classical solution of the corresponding HJB  equation, which in turn leads to  an optimal restricted payment policy known as the threshold strategy. In the case of unrestricted payment scheme,  by solving the associated integro-differential quasi-variational inequality, we obtain the value function as well as an optimal unrestricted dividend payment scheme known as the barrier strategy. When claim sizes are exponentially distributed,  we provide easily verifiable conditions under which the threshold and barrier strategies are optimal restricted and unrestricted dividend payment policies, respectively. The main results are illustrated with several examples, including a new example concerning regressive growth rates.



\end{abstract}

 \setlength{\baselineskip}{0.201in}
\section{Introduction}
\label{sect-Introduction}
 The dividend problem in classical insurance risk models was originated in \cite{DeFi}, and drew revived interests in  
 recent literature focusing on optimization of dividend payment strategies. The optimality is often considered to be a strategy which maximizes the expected present value of dividends received by the shareholders.
 \cite{Jeanblanc1995} and  \cite{Asmussen1997} investigated in diffusion models the dividend problems where the dividends are permitted to be paid out up to a maximal constant rate or a ceiling. We shall refer to such a type of dividend problem as restricted payment scheme. It was shown in their papers that the dividends
should be paid out at the maximal admissible rate as soon as the
surplus exceeds a certain threshold. Interestingly, it turns out
that such a threshold strategy is the optimal restricted payment scheme in a variety of other risk models. For example,
 \cite{Gerber2006}  discussed the threshold strategy in the compound
Poisson model and solved the problem explicitly when the claim size
is exponentially distributed. \cite{Fang2007} studied a similar problem in the compound Poisson risk model with constant
interest and showed the optimal dividend strategy is a threshold
strategy for the case of an exponential claim distribution.
See also \cite{Asmussen-H-T-00},  \cite{Bai-Paulsen-10},
  \cite{Choulli}, \cite{HunPau}
   \cite{Schmidli-02}, and references therein
 for some important developments on optimal dividend policies in the setting of controlled diffusions.

On the other hand, there was also a significant amount of literature in which no such restriction of maximal rate is imposed on dividend payment strategies.
We shall refer to this type of dividend strategy as the unrestricted payment scheme.
% Unrestricted payment
Such schemes are motivated by the fact that
 dividends are not usually paid out in a continuous fashion in practice.
For instance, insurance companies may distribute dividends on discrete time intervals, in theory % allowing for 
resulting in unbounded payment rate. In such a scenario,
 the surplus level changes drastically on a dividend payday.
  In other words, abrupt or discontinuous changes occur due to
  ``singular'' dividend distribution policy.
  This gives rise to a singular stochastic control problem.
% It was shown universally in various risk models that the optimal strategy
% is to pay nothing until the surplus reaches a certain level at which
% point dividends should be paid so as to keep the surplus from exceeding
% this level. Such a strategy is also known as the barrier strategy in
% ruin literature or the band strategy in stochastic controls.
Such problems are studied in   \cite{Choulli,Paulsen-08,Shreve-84}, % \cite{Paulsen-07,Paulsen-08},  %  \cite{Paulsen-97},
 and the references therein
when the surplus dynamics is modeled by a controlled diffusion.
But to the best of our knowledge,
related work in the setting of piecewise-deterministic compound Poisson risk model is relatively scarce. The most notable includes
 \cite[Section 2.4]{schmidli2008} and \cite{Albrecher-Thon-08}, which
investigate the optimal unrestricted dividend payment problem when the surplus process follows a classical Cramer-Lundberg risk model without or with the force of interest, respectively.

As pointed out in \cite{Cai-Feng},
the classical Cram\'er-Lundberg risk model and the compound Poisson risk model with interest, and the compound Poisson risk model with
absolute ruin are all special cases of piecewise-deterministic compound Poisson
(PDCP) risk model.
One naturally asks
 whether there exist unifying optimal solutions to both
 dividend payment schemes in PDCP  risk models. 
 % If so, can we confirm in general that 
 Moreover, can we find the most general conditions under which the threshold strategy is the optimal restricted dividend policy whereas the barrier strategy is the optimal unrestricted dividend policy?
%We provide affirmative solutions to both of these
% questions in this paper under certain conditions.
 We formulate and solve the problems within the framework
of stochastic control theory in the specific setting of PDCP risk model.
%Roughly, the idea is to pay out the dividend
%at a dynamic rate in such a way that a certain reward function is optimized.
 Compared with the aforementioned work in the setup of controlled diffusions,
 the associated Hamilton-Jacobi-Bellman (HJB) equation in our work contains a non-local term
(the integral term with respect to the claim size distribution), resulting in
substantial difficulty and technicality in the analysis.
% Nevertheless, we use renewal type arguments to overcome this difficulty
% and establish the HJB equation.
  % and prove that the value function in the restricted
  % payment scheme is a classical solution to
  % the HJB equation \eqref{HJB}.

The contribution and novelty of this work  also arise from several different aspects.
\begin{enumerate}

    \item A salient feature of our model is the generality of pure jump models in which both restricted 
    and unrestricted payment schemes are presented and directly compared.
    Although special cases of the  % piecewise-deterministic compound Poisson
    PDCP risk model have been treated
      in the literature, this paper extends enormously the spectrum of risk models which exhibit common optimality.

    \item  We obtain general optimal solutions in the case of exponential claim size distribution.      Moreover, we provide sufficient conditions to guarantee the optimality of the threshold and barrier strategies for the restricted and unrestricted dividend payment schemes in a general PDCP risk model.  
 To the best of our knowledge,    these conditions were unknown previously  in the literature.  
 Note that the analysis of solutions
 in the general case (Theorems \ref{finalthr} and \ref{finalbar}) in this paper is entirely based on qualitative study of ordinary differential equation (ODE) and integro-differential equation (IDE).
	
 % Second, the generality of solution formulation by which  explicit
 % solutions can be obtained in the case of exponential claim size distribution. It turns out

    \item It is also worth mentioning that the solution methods presented in this paper are shown with examples to be more efficient alternatives to the known methods in the existing literature. For example, we propose simple procedures (Theorems \ref{resthm} and \ref{thm-soln-qvi})  to identify the optimal threshold and barrier levels, rather than the optimization procedure on value functions which would have to be first determined explicitly. 

 %   \item Finally, we point out the solution methods in this paper are particularly suitable for implementation, taking advantage of modern symbolic calculation systems. Many previous works on similar topics rely on solving systems of equations satisfying miscellaneous boundary conditions, which often require intellectual input for simplification. In contrast, the optimal value functions in Theorems \ref{resthm} and \ref{thm-soln-qvi} can be derived systematically within merely a few lines of automated procedure.

\end{enumerate}

The rest of the paper is organized as follows. After recalling the notion of a PDCP model, we formulate the   optimality of dividend strategies   as a
stochastic control problem in Section \ref{sect-model}. We consider in Section \ref{sect-restricted} the restricted dividend
payment schemes. Some properties of the value function are derived
 and the value function is shown to be a classical solution to the HJB equation \eqref{HJB}.
 In Section \ref{sect-singular}, we formulate the optimal unrestricted payment scheme
 as a singular stochastic control problem and establish a
 verification theorem of the quasi-variational inequality \eqref{quasi-variational-ineq}.
 % is established.  Furthermore,
% under some general conditions, we provide an explicit procedure to obtain
% an optimal dividend barrier and the corresponding optimal value function.
When the claims are exponentially distributed, with complete generality of the PDCP risk model, we provide easily verifiable sufficient conditions
for the  optimality of threshold and barrier  dividend payment schemes and
  obtain explicit solutions
for both the restricted and unrestricted dividend payment schemes in Sections \ref{sect-restricted} and
\ref{sect-singular}, respectively.    %is devoted to
 Three % explicit
 examples are % produced
 provided  for illustrative purpose in Section \ref{examples}.
 Finally, the paper is concluded with several remarks in Section \ref{sect-conclusion}. Two technical results on the qualitative analysis of the solution to an ODE and several proofs are  placed in the Appendix. 
 
A preliminary version of the paper was announced in \cite{pdcp-cdc-2012}  without proofs. 
In addition, the current version contains new results on sufficient conditions for the optimality of the threshold and barrier strategies (Theorems  \ref{finalthr} and \ref{finalbar}).  

To facilitate later presentations, we introduce some notations here.
We use $I_A$ to denote the indicator function of a set $A$. When $a,b \in \rr$, $a\wedge b:=\min\{a,b\}$ and $a^+:=\max\{a, 0\}$.
Throughout the paper, we use the notations $\xi(t)$ and $\xi_t$ interchangeably.
 % A function $\xi$ from $[0,\infty)$ to some Polish space $E$ is {\em c\`adl\`ag} if it is right continuous and has left limits in $E$. 
 % If $E=\rr$ and $\xi $ is c\`adl\`ag, then $\Delta \xi(t)= \xi(t)-\xi(t-)$ for $t> 0$ and the convention $\Delta \xi(0)= \xi(0)$ is used.
 % The continuous part of $\xi $ is denoted by $\xi^c(t):=\xi(t)- \sum_{0\le s \le t} \Delta \xi(s)$.
As usual, $\sup  \emptyset =-\infty$ and $\inf \emptyset = + \infty$.

\subsection{Problem Formulation}\label{sect-model}

To give a rigorous formulation of the optimization
problem, we start with a filtered probability space $ \{
\Omega,\mathcal {F}, {\{\mathcal {F}_t }\}_{t\geq0},\mathrm{P} \}$ satisfying the usual condition.
We assume that the surplus level  is modeled by a
piecewise-deterministic compound Poisson process. Note that the jump points represent the arrivals of insurance claims and the downward jumps are determined by claim sizes.

%%%%%% 
Suppose the surplus level $X(t)$ of an insurance company at time $t \ge 0$ is modeled by a  piecewise-deterministic compound Poisson process
\begin{equation}\label{uncontrolled-surplus}
   X(t)  =x+ \int_0^t g(X(s))ds -\sum_{i=1}^{N(t)} Y_i,  \ \ t \ge 0, 
\end{equation} where $x\ge 0$ is the initial surplus, $N=\set{N(t),t\ge 0}$ is a Poisson process with rate $\lambda>0$, $Y_1, Y_2, \dots$ are independent and identically distributed nonnegative random variables, and $g$  is a Lipschitz continuous function, taking values in $(0,\infty)$, and satisfies the linear growth condition. Denote the common distribution function of $Y_1, Y_2, \dots$ by $Q$.

Denote by $0=T_0 < T_1 < T_2 < \cdots $ the sequence of jump points of the process $X$, then $\Delta X({T_k}):= X({T_k})-X({T_k-})= Y_k$ for $k=1,2,\dots$. Also, the surplus process $X$ between any two consecutive jumps is deterministic and given by 
$ X (t)= \phi_{X({T_k})}(t-T_k),$   $   t\in [T_k, T_{k+1}),$  $   k=0, 1, 2, \dots, $
    with $\phi_z(t)$  determined by
    $ \df \phi_z(t) = g(\phi_z(t))\df t, t > 0$
and $\phi_z(0)=z$. 

% \begin{defn}\label{PDCP}
% A piecewise-deterministic compound Poisson (PDCP) process is a
% real-valued stochastic  process  $ X=\{ X({t}),0 \leq t  < \infty   \}  $,
% defined on a given  probability space  $   \{ \Omega,\mathcal {F}, \set{\F_t}_{t\ge 0},
% \mathrm{P} \}$, satisfying the following properties:
% \begin{itemize}
% \item[(i)] $X(0)=x\ge 0 $,
% \item[(ii)] Let $0=T_0 < T_1 < T_2 < \cdots $ denote a sequence of jump points of the process $X$. 
% Assume  that $T_{i+1}-T_i$ has exponential distribution with mean $1/\la>0$ for every $i=0,1,\dots$
% % Then the adapted counting process defined by $N(t)=\sum_{i=1}^\infty I_{\set{T_i\le t}}$
% %  follows a homogeneous Poisson process with intensity rate $\la$,
% \item[(iii)] The jump sizes $Y_k=\Delta X({T_k})= X({T_k})-X({T_k-})$ for $k=1,2,\dots$ are independent 
% and identically distributed nonnegative random variables with common distribution function
% $Q(y)=1-\bar Q(y)= \pr\set{Y_1\le y},  0 \le y < \infty$,
% \item[(iv)] The process between any two consecutive jumps is deterministic and given by
%    $ X_t= \phi_{X({T_k})}(t-T_k), \ \ t\in [T_k, T_{k+1}), \ \ k=0, 1, 2, \dots, $
%    where $\phi_z(t)$ is determined by
%    $ \df \phi_z(t) = g(\phi_z(t))\df t, t > 0, $
 %   satisfying $\phi_z(0)=z$ and $\lim_{t\to \infty}\phi_z(t) =L \in [-\infty, \infty]$. 
% $ The function $g:  B\to (0,\infty)$, satisfies the linear growth condition and is
%   Lipschitz continuous on %each subinterval on a finite partition of
%     its domain $B$.
% \end{itemize}
% \end{defn}

We give the corresponding expressions for $ g(x)$ and
$ \phi_{x } (t)  $ for three special cases in which optimal dividend policies will be developed later as examples.
\begin{itemize}
    \item (Cram\'er-Lundberg model) The deterministic growth in surplus
between any two consecutive claims is defined by the influx of premium at a constant rate $c$ per time unit, i.e., $g(x)=c,~x \geq 0$.
Hence, $\phi_{x } (t) =x+ct, ~t \geq 0 $.

    \item (Constant interest model) All positive surplus earns interest at the constant rate $\rho > 0$ per time unit, i.e., $g(x)=
\rho x +c, ~x \geq 0 $. Hence, $\phi_{x } (t) = ( x + c/ \rho) e
^{\rho t }{ - c/\rho} , ~t \geq 0  $.

	\item (Regressive growth model) The surplus growth rate tends to regress to a mean premium rate $c>0$ according to $\df g(x)=b\big(c-g(x)\big)\df x$ with $g(0)=a+c$ and  $a,b>0.$ Hence, $g(x)=ae^{-bx}+c$ and $\phi_x(t)=(-1/b)\{\ln c-\ln[(ae^{-bx}+c)e^{b(x+ct)}-a]\}$ for $t\ge 0$.

%    \item (Classical model with absolute ruin) All positive surplus is invested with constant interest return at the lending rate $\rho>0$ per time unit whereas the negative surplus are kept with constant interest payment at the borrowing rate $r>0$ per time unit. Ruin occurs only the incoming premium is no longer large enough to sustain the outgoing interest payment. In this case, $t_0=(1/r)\log(c/(rx+c))$ and
%        \[g(x)=\left\{  \begin{array}{ll} \rho x+c, & x \ge 0;\\ rx+c, & -c/r<x<0,\end{array} \right. \phi_x(t)=\left\{  \begin{array}{ll} (x+c/\rho)e^{\rho t}-c/\rho, & x \ge 0, t \ge 0;\\ c/\rho(e^{\rho(t-t_0)}-1) , & -c/r<x<0, t\ge t_0\\ (x+c/r) e^{rt}-c/r,& -c/r<x<0, 0 \le t< t_0 \end{array} \right.\]
\end{itemize}
More PDCP models can be found in \cite{Cai-F-W-09, Cai-Feng, Cai-Feng2}, and \cite{Alb-Har-2007}.

We now  enrich the model
by considering   dividend payout. Denote by $D(t)$ the aggregate dividends by time $t$. Assume that $D=\set{D(t), t\ge 0}$ is c\`adl\`ag, nondecreasing, and $\F_t$-adapted with $D({0-})=0$.
 Moreover, we require that at any time $t$, the dividend payment should not exceed the current surplus level, i.e., 
 $\Delta D(t):=D(t)-D(t-)\leq X^D(t-).$
Any dividend payment scheme $D=\set{D(t), t\ge 0}$ satisfying the above conditions
is called an {\em admissible control} and the collection of all admissible controls is denoted by $\Pi$. %$=\set{D: D \text{ is admissible}}$.
The dynamics of the controlled surplus process under the admissible control $D$ is
 \beq\label{controlled-surplus}
X^D(t)= x+ \int_0^t g(X^D(s))ds -\sum_{i=1}^{N(t)} Y_i -D(t), 
\eeq where $x\ge 0$ is the initial surplus.
The time of ruin is denoted by $$\tau=\tau(x,D):=\inf\set{t\ge 0: X^D(t) < 0}.$$
The %performance functional or the
 expected present value (EPV) of dividends up to ruin is defined as
\beq\label{performance-defn}
J(x,D)= \ex_{x} \int_0^\tau e^{-\dl t} \df D(t),
\eeq
where $\dl >0$ is the force of interest. The objective is to find an admissible control $D^*\in \Pi$ that maximizes the EPV. That is, we seek 
\beq\label{value-defn} V(x):= \sup_{D \in \Pi} \set{J(x,D)} = J(x,D^*). \eeq
Note that $V(x)=0$ for all $x< 0$.
Depending on the parameters of the model, $V$ can be $\infty$. In the rest of the paper, to work with a well-formulated maximization problem, we assume that $V(x)< \infty$ for all $x\ge 0$. See Section \ref{sect-singular} for  a sufficient condition for the finiteness of the value function.

\section{Restricted Payment Scheme}\label{sect-restricted}

We first consider problem \eqref{value-defn} for the case when the dividend payment scheme $D=D_R$ is absolutely continuous \wrt time. That is, there exists some $u(t), t\ge 0$ such that
$ D_R(t) = \int_{0}^t u(s) \df s.$
Moreover, we assume that $u(t)$ is $\F_t$-adapted and that there exists some
positive constant $u_0$ % with $u_0 < g(x)$ for all $x \ge 0$ 
such that
$ 0\le u(t) \le u_0< \inf \set{g(x),x \ge 0}$,  for all $ t \ge 0. $
Denote the collection of all such dividend payment schemes by $\Pi_R$.
The EPV corresponding to the initial surplus $x \ge 0 $ under the dividend payment policy $D_R=\set{D_R(t), t\ge 0}$ is given by
\beq\label{epv-restricted}
J(x,D_R) = \ex_x \int_0^\tau e^{-\dl t} \df D_R(t) =\ex_x \int_0^\tau e^{-\dl t} u(t) \df t.
\eeq
The goal is to find an admissible policy $D^{*}_R \in \Pi_R$ such that \beq\label{restricted-V-defn}
  V_R(x):=\sup_{D_R \in \Pi_R} J(x,D_R)=J(x,D^{*}_{R}).\eeq
Apparently, we have $V_R(x) \le V(x)$ for all $x\ge 0$, where $V(x) $ is the value function in \eqref{value-defn}.

 \subsection{The HJB Equation  and  Optimal Strategies}

% We need  the following {\em dynamic programming principle}.
% For any $ \mathcal {F}_{t}$-stopping time $\theta$,
%\begin{eqnarray}\label{DPP}
% V_R (x)=\sup_{u\cd\in \Pi_R}
%\mathrm{E}_{x}\left[\int_0^{ \tau \wedge \theta } e^{-\delta s}
% u(s) \df s +  e^{-\delta (  \theta \wedge  \tau)}V_R(X (\theta \wedge \tau  ) )   \right],\ \ x\geq 0
%\end{eqnarray}
 % The dynamic programming principle is well-known in the  diffusion case (see, for example,  \cite{FlemingS} and  \cite{Krylov}) 
 % and is also proved for jump diffusions in   \cite{Ishikawa-04}. 
% Various papers such as  \cite{Davis} and \cite{Azcue2005} employed the dynamic programming principle 
% in the setting of piecewise-deterministic Markov processes, which also enables us to derive the HJB equation in this section. 

 We first derive some elementary properties of the value function \eqref{restricted-V-defn}, which  will
 enable us to establish the HJB equation in \thmref{thm-HJB}. % To % showcase the more concrete results, 
The proofs can be found in the Appendix.

 \begin{lem}\label{lem-bounds}
The function $V_R(x)$  is bounded by $u_{0} / \delta $,
increasing, and Lipschitz continuous on $[0,\infty)$, and therefore
absolutely continuous, and converges to $ u_{0}/ \delta   $ as $x
\rightarrow  \infty$.
 \end{lem}

\begin{thm}\label{thm-HJB} The function
  $V_R(x)$ is differentiable and
fulfills the HJB equation
\begin{equation}\label{HJB}
 \sup_{0\leq u  \leq u_0 }\bigg
\{[g(x)-u ]V_R'(x)-(\lambda+\delta)V_R(x)   +\lambda\int^{x}_{0}V_R(x-y)\df Q(y)+u \bigg\}=0, \quad x \ge 0.
\end{equation} 
Moreover, the  strategy $D^{*}_{R}=\{\int^t_0 u^{*}_{R}(s)\df s, s\geq 0\}$ with
  \begin{equation}\label{eq:u_R}
u^{*}_{R}(t)= \begin{cases}
0, & \text{ if }\  V_R'(X^{*}_{R}(t))>1,\\
u_{0},& \text{ if }\  V_R'(X^{*}_{R}(t))\leq 1,
\end{cases}
\end{equation} is optimal in the sense that $J(x,D^{*}_{R})=V_R(x),$
where $X^{*}_{R}(t)$  
 is  the
corresponding surplus process under the strategy $D^{*}_{R}$.
\end{thm}

%\begin{proof} The proof of this theorem utilizes Lemma \ref{lem-bounds}, some renewal type argument,  and the dynamic %programming principle: 
%\begin{eqnarray}\label{DPP}
%V_R (x)=\sup_{u\cd\in \Pi_R}
%\mathrm{E}_{x}\left[\int_0^{ \tau \wedge \theta } e^{-\delta s}
% u(s) \df s +  e^{-\delta (  \theta \wedge  \tau)}V_R(X (\theta \wedge \tau  ) )   \right],\ \ x\geq 0,
%\end{eqnarray}  where $\theta$ is an $ \mathcal {F}_{t}$-stopping time.
% For brevity we shall omit the details here. Similar argument in the case of  classical Cram\'er-Lundberg model can be found in    %\cite[Theorem 2.32]{schmidli2008}. 
%  \end{proof}
  
%The proof of this theorem requires lengthy technical details. For the clarity of the presentation, we leave the proof in a supplement posted on the authors' website.

\subsection{Exponential Claims}\label{sect-exponential1}
We consider a simple yet thought-provoking case where an explicit general solution to the HJB equation \eqref{HJB}
  and an optimal dividend payment policy can be obtained. Assume that the common claim size distribution is given by
$  Q(y) = 1- e^{-\al y},   y \ge 0, $   for some  $ \alpha>0. $ 

\begin{lem}\label{sec34-lem1} The IDE
\beq g(x)\varphi'(x)-(\lambda+\delta)\varphi(x)+\lambda \int^x_0 \varphi(x-y) \alpha e^{-\alpha y} \df y=0, \qquad x>0,\label{intdiff}\eeq
has a positive and strictly increasing solution $\psi_{1}$.
\end{lem}
\begin{proof}
It is well known that the IDE \eqref{intdiff} has a unique solution $\psi_{1}: [0,\infty) \to \rr$ determined up to a multiplicative constant.
 Without loss of generality, we choose $\psi_{1}(0)=g(0) > 0$. Then \eqref{intdiff} implies that $\psi'_{1}(0)= \lambda + \delta > 0$.
 Now the desired assertion will follow if we can prove that $\psi_{1}'(x) > 0$ for all $x\ge 0$.
Suppose this was not the case, then there would exist some $x_{1}>0 $ such that $\psi_{1}'(x) > 0 $ for all $x\in [0,x_{1})$ but $\psi'_{1}(x_{1}) = 0$. Then by virtue of \eqref{intdiff},
\bea 0\ad = g(x_{1})\psi_1'(x_1) - (\lambda + \delta )\psi_1(x_1) + \lambda \int_0^{x_1} \psi_1(x_1-y) \alpha e^{-\alpha y} \df y \\
           \ad      < - (\lambda + \delta )\psi_1(x_1) + \lambda \int_0^{x_1} \psi_1(x_1) \alpha e^{-\alpha y} \df y  <      -\delta \psi_1(x_1)  < 0.
    \eea
   This is a contradiction and therefore we must have $\psi_{1}'(x) > 0$ for all $x\ge 0$.
\end{proof}

% In addition, we make the following assumption.
%\nd {\bf Hypothesis A} The integro-differential equation
%\beq g(x)\varphi'(x)-(\lambda+\delta)\varphi(x)+\lambda \int^x_0 \varphi(x-y) \alpha e^{-\alpha y} \df y=0, \qquad x>0,\label{intdiff}\eeq has a strictly increasing
% and concave
% solution $\psi_1(x)$ and the  differential equation
%\beq [g(x)-u_0]\varphi''(x){+}[\alpha g(x)-\alpha u_0+ g'(x)-(\lambda+\delta)]\vphi'(x)-\alpha \delta \vphi(x)=0, \qquad x>0,\label{intdiff2}\eeq
%has a bounded concave solution $\psi_2 (x)$.

\begin{thm} \label{resthm}
Let $\psi_1$ be  a positive and strictly increasing solution  %determined up to a constant factor 
to the IDE
\eqref{intdiff}. 
Suppose  that $\psi_2$ is a strictly increasing and concave solution %determined up to a constant factor 
to the  ODE
\beq [g(x)-u_0]\varphi''(x){+}[\alpha g(x)-\alpha u_0+ g'(x)-(\lambda+\delta)]\vphi'(x)-\alpha \delta \vphi(x)=0, \qquad x>0,\label{intdiff2}\eeq
and  that there exists a
%unique
number $d>0$ such that
$\psi_1$ is concave on $(0, d)$ and
%$\psi_2'(d) > 0$,
\beq  \frac{\psi_1(d)}{\psi_1'(d)}-\frac{\psi_2(d)}{\psi_2'(d)}=\frac{u_0}{\delta}.
\label{condn}\eeq
Then the %solution to \eqref{HJB}
value function $V_R(x)$ is given by
\begin{equation}\label{solnhjb} V_R (x)  =\begin{cases} \displaystyle \frac{\psi_1(x)}{\psi'_1(d)}, & \text{ if }\ 0 \le x <  d,\\
\displaystyle \frac{u_0}{\delta} +  \frac{\psi_2(x)}{\psi'_2(d)}, & \text{ if }\ x \ge  d.\end{cases} \end{equation}
 Moreover, the optimal dividend payment policy is the {\em threshold strategy}
 \beq\label{optimal-restricted-policy} u^*(t)=u_0I_{\set{ X^*(t) \ge  d}}, \eeq
 where $X^*$ is the corresponding controlled surplus process.
\end{thm}

The interpretation of such an optimal strategy
 is as follows.
First, a
threshold $d$ is determined so that dividend payments start immediately when the threshold is attained. Second, as long as the surplus process remains above the threshold, dividends are paid out continuously at the maximal rate $u_0$ per time unit. No dividend payment is  allowed when the surplus drops below the threshold.

\begin{proof}
% Without loss of generality, by virtue of Lemma \ref{sec34-lem1}, we can assume $\psi_1$ is a positive and increasing function.
% the existence of which is shown in Lemma \ref{sec34-lem1}.
Denote by $\Psi(x)$ the function defined on the right-hand side of \eqref{solnhjb}. Note that $\Psi $ is continuously differentiable with $\Psi'(d)=1$.
 Since both $\psi_1$ and $\psi_2$ are concave functions and $\psi'_1(d)>0, \psi'_2(d)>0$, we must have $\Psi'(x)>1$ for $0\le x<d$ and $\Psi'(x)<1$ for all $x>d.$
Hence by virtue of Theorem  \ref{thm-HJB}, %and \ref{optimalvalue},
 it only remains to show that $\Psi$ satisfies the % integro-differential
HJB equation \eqref{HJB}.

%First of all,
It is clear by definition that
\[ g(x)\Psi'(x)-(\lambda+\delta)\Psi(x)+\lambda \int^x_0 \Psi(x-y) \alpha e^{-\alpha y} \df y=0, \qquad 0\le x <d.\]
Therefore $\Psi $ solves the HJB equation \eqref{HJB} if
  %Second, we want to
  we can show that
\begin{equation} [g(x)-u_0]\Psi'(x)-(\lambda+\delta)\Psi(x)+\lambda \int^x_0 \Psi(x-y) \alpha e^{-\alpha y} \df y+u_0=0, \qquad x\ge d.\label{eqn2}\end{equation}
To this end, we define
$ h(x)=\alpha e^{-\alpha x}\int^x_d e^{\alpha y} \psi_2(y)\df y,$  % \quad
for $  x\ge d. $
It is straightforward to verify that
 $\alpha \psi_2(x)=h'(x)+\alpha h(x).$
 % Denote the left-hand side of \eqref{intdiff2} by LHS. Then 
 Consequently, \eqref{intdiff2} can be rewritten as 
\bea
0 \ad =[g(x)-u_0]\psi''_2(x){+}[\alpha g(x)-\alpha u_0+ g'(x)-(\lambda+\delta)]\psi'_2(x)-\alpha (\lambda+\delta) \psi_2(x) +\lambda h'(x)+\lambda \alpha h(x)\\
\ad = [g(x)-u_0]\psi''_2(x)+g'(x)\psi'_2(x)-(\lambda+\delta)\psi'_2(x)+\lambda h'(x)+ \\ \aad \hfill  \alpha\{ [g(x)-u_0]\psi'_2(x)- (\lambda+\delta) \psi_2(x)+\lambda h(x)\}.
\eea
Denote the LHS of \eqref{eqn2} by $H(x)$. The equation above shows that $H'(x)+\alpha H(x)=0$ for all $x\ge d$. Letting $x=d$ in \eqref{intdiff} and using $\Psi'(d)=1$, we see that $H(d)=0$. Thus, $H(x)=0$ for all $x\ge d$ and the claim \eqref{eqn2} is proved.

Finally the verification of optimality for the strategy defined in \eqref{optimal-restricted-policy} is straightforward and we shall omit the details here.   
\end{proof}

Theorem \ref{resthm} establishes the optimality of the threshold strategy and provides an easy procedure to identify the threshold level $d$. These results are based on the existence of solutions to the IDE \eqref{intdiff} and the ODE \eqref{intdiff2} with  certain properties.
In particular, the smooth pasting condition \eqref{condn} must hold. One may naturally ask under what conditions these solutions and the threshold level $d$  exist.
 The following  theorem gives a minimal set of easily verifiable conditions in the PDCP model. % with complete generality.

\begin{thm} \label{finalthr}
Suppose that $g\in C^2([0,\infty))$ and satisfies
\begin{align} \label{eq:g-0} \alpha \lambda g(0)+(g'(0)-\lambda-\delta)(\lambda+\delta)>0;\\
 \label{eq:sup-g-x} \displaystyle \sup_{x \ge 0} \left\{g''(x)+\alpha g'(x)-\alpha \delta \right\}<0.\end{align}
 Then \eqref{intdiff2} admits a solution $\psi_2$  that is negative, strictly increasing and concave.
 
 \begin{itemize}\item[(i)] Furthermore, if
\begin{equation}\label{unique}  \frac{\psi_2(0)}{\psi'_2(0)} > \frac{g(0)}{\lambda+\delta}-\frac{u_0}{\delta}, \end{equation}
then there exists a unique $d>0$ such that equation \eqref{condn} holds and consequently the optimal restricted payment scheme is \eqref{optimal-restricted-policy} and the value function is given by \eqref{solnhjb};

\item[(ii)] Otherwise, if
\begin{equation}\label{Kcond}\frac{\psi_2(0)}{\psi'_2(0)}\le \frac{g(0)}{\lambda+\delta}-\frac{u_0}{\delta},\end{equation} then the optimal restricted payment scheme is $D^\ast_R=\{\int^t_0 u^\ast_R(s) \df s, s \ge0\}$ with $u^\ast(t)=u_0$ for all $t\ge 0$ and the value function is given by
\beq\label{V_R-K}V_R(x)=\frac{\psi_2(x)}{K}+\frac{u_0}{\delta},\qquad x \ge 0,\eeq where
\[K=\frac{[g(0)-u_0]\psi'_2(0)-(\lambda+\delta)\psi_2(0)}{(\lambda/\delta)u_0}.\]
\end{itemize}
\end{thm}

\begin{proof}
 %We now prove the assumptions of Theorem \ref{resthm} are satisfied.
Take $f(x)=g(x)-u_0$, $ h(x)=g'(x)+\alpha g(x)-\alpha u_0-(\lambda+\delta)$, and $ k=\alpha \delta$ in \eqref{sode}.
 % Thus $\psi_2(x)$ satisfies $f(x)\psi''_2(x)+h(x) \psi'_2(x)-k \psi_2(x)=0,$ for $x>0$.
  Note that $f(x) > 0$ and by virtue of \eqref{eq:sup-g-x},  $h'(x) = g{''}(x) + \alpha g'(x) \le m < k$ for all $x \ge 0$, where $m=  \sup_{x\ge 0} \set{g{''}(x) + \alpha g'(x)} $.
  %As proved previously, there exists $0<m<k$ such that $h'(x) \le m<k$ for all $x \ge 0$.
  It then  follows from Lemma \ref{odecond} that there exists a negative, strictly increasing and concave solution $\psi_2(x)$ to \eqref{intdiff2}  on $[0, \infty).$

(i) In this case, by virtue of Theorem \ref{resthm}, it remains to prove that there is a $d>0$ such that  $\psi_1(x)$ is concave on $(0,d)$ and that \eqref{condn} is satisfied, where $\psi_1$ is a solution to the IDE \eqref{intdiff}. As argued in \lemref{sec34-lem1}, we can take $\psi_1(0)=g(0)$ and $\psi'_1(0)=\lambda + \delta$.
Next, applying the operator $(\df/\df x+\alpha)$ to \eqref{intdiff} gives the ODE
\beq g(x) \psi''_1(x)+(g'(x)+\alpha g(x)-\lambda-\delta) \psi'_1(x)-\alpha \delta \psi_1(x)=0.\label{eq-psiode}\eeq
 In particular, due to \eqref{eq:g-0}, letting $x=0$ in \eqref{eq-psiode} yields that
\bed \begin{aligned} g(0)\psi''_1(0)&=\alpha \delta g(0)-(g'(0)+\alpha g(0)-\lambda-\delta)(\lambda+\delta)
=-\alpha \lambda g(0)-(g'(0)-\lambda-\delta)(\lambda+\delta)<0. \end{aligned} \eed
 Therefore, $\psi''_1(0)<0.$
With $f(x)=g(x)$ and $h(x)=g'(x)+\alpha g(x)-\lambda-\delta$, as argued before, $h'(x) \le m  < \alpha \delta. $
Now it follows from Lemma \ref{idecond} that  there exists a $b>0$ such that $\psi_1''(b)=0$, $\psi_1''(x)<0$  for $x \in [0,b)$ and  $\psi_1''(x)>0$  for $x \in (b,\infty)$.

Using \eqref{intdiff2} and \eqref{eq-psiode}, we obtain
$$(g(x)-u_0)\frac{\psi_2''(x)}{\psi_2'(x)} = g(x) \frac{\psi_1''(x)}{\psi_1'(x)}
+ \alpha \delta \left( \frac{\psi_2(x)}{\psi'_2(x)} + \frac{u_0}{\delta} - \frac{\psi_1(x)}{\psi_1'(x)}\right), \ \ \forall x \ge 0.$$
 In particular, noting $\psi_1''(b)=0$, the above equation yields
 $$\frac{\psi_2(b)}{\psi'_2(b)}  + \frac{u_0}{\delta} - \frac{\psi_1(b)}{\psi_1'(b)}  = \frac{1}{\alpha\delta}(g(b)-u_0)\frac{\psi_2''(b)}{\psi_2'(b)}  < 0,$$
where the last inequality follows from the fact that $\psi_2$ is increasing and concave.
On the other hand, it follows from \eqref{unique}
 that
 $$\frac{\psi_2(0)}{\psi'_2(0)}  + \frac{u_0}{\delta} - \frac{\psi_1(0)}{\psi_1'(0)}
 =  \frac{\psi_2(0)}{\psi'_2(0)}  + \frac{u_0}{\delta} - \frac{g(0)}{\lambda + \delta} > 0. $$
 By the intermediate value theorem, the solution $d$ of \eqref{condn} exists and $d  \in (0, b)$.

 The uniqueness of $d\in (0, b) $ follows immediately from the mean value theorem.
 Suppose on the contrary that there were $d_1<d_2 \in (0,b)$ both satisfying \eqref{condn}. Then there would exist some $\xi \in (d_1,d_2)$ with $$\frac{\df }{\df x} \(\frac{\psi_1(\xi)}{\psi'_1(\xi)}\) = \frac{\df }{\df x} \(\frac{\psi_2(\xi)}{\psi'_2(\xi)}+\frac{u_0}{\delta} \).  $$
 But this is impossible, since
 % On one hand, since $\psi_1(x)>0$ for all $x\ge 0$ and $\psi''_1(x)<0$ for $0 \le x <b,$  for $0\le x <b$,
\[\frac{\df }{\df x} \left( \frac{\psi_1(x)}{\psi'_1(x)} \right)=1-\frac{\psi''_1(x) \psi_1(x)}{(\psi'_1(x))^2}>1, \ \ \forall x \in (0,b), \]
% On the other hand, since $\psi_2(x)<0$ and $\psi''_2(x)>0$ for all $x \ge 0$, then
and
\[\frac{\df }{\df x} \left( \frac{\psi_2(x)}{\psi'_2(x)}+\frac{u_0}{\delta} \right)=1-\frac{\psi''_2(x) \psi_2(x)}{(\psi'_2(x))^2}<1 \ \ \forall x \in (0,b). \]

(ii) Note that the condition \eqref{Kcond} is equivalent to $K \ge \psi'_2(0).$ Since $\psi''_2(x)<0$ for all $x \ge 0$,  we have  $\psi'_2(x)< \psi'_2(0) \le K$ and hence $\psi'_2(x)/K \le 1$ for all $x \ge 0$. Moreover, as  $\psi'_2(x)/K$ satisfies the boundary condition
\[[g(0)-u_0]\frac{\psi'_2(0)}{K}-(\lambda+\delta)\left(\frac{\psi_2(0)}{K}+\frac{u_0}{\delta} \right)+u_0=0,\] using an argument similar to that in the proof of \thmref{resthm}, it is easy to verify  that  $V_R(x)$ given in \eqref{V_R-K} is indeed a solution to the IDE \eqref{eqn2} and hence a solution to the HJB equation \eqref{HJB}.
\end{proof}

% Last, we prove the existence and uniqueness of $d \ge 0$ under the condition \eqref{unique}. Note that $d=0$ when
% \[\frac{\psi_1(0)}{\psi'_1(0)}= \frac{g(0)}{\lambda+\delta}= \frac{\psi_2(0)}{\psi'_2(0)} +\frac{u_0}{\delta}.\]

\begin{rem}
With the complete generality of the PDCP model, we have shown that under assumptions \eqref{eq:g-0} and \eqref{eq:sup-g-x},
 the optimal restricted dividend payment scheme is always the threshold strategy, which is to pay dividends at the maximal rate as long as the surplus is above a threshold level. In the case  % that % $\psi_2(0)/\psi'_2(0)> g(0)/(\lambda+\delta)-u_0/\delta$, 
	of \eqref{unique}, the threshold level is a unique positive number. Otherwise, if \eqref{Kcond} holds, the threshold level is set to be zero.
\end{rem}

\section{Unrestricted Payment Scheme}\label{sect-singular}
In   Section \ref{sect-restricted}, we   considered the case where the dividend payment rate is bounded.
Consequently, the surplus level changes continuously in time $t$ in response
to the dividend payment policy.
However, in many applications,
the boundedness of the dividend payment rate
seems rather restrictive. For instance, insurance companies are more likely to distribute the dividend at discrete time points rather than with a continuous stream of dividend payments. Thus we remove the restriction on the maximal dividend rate and  consider the (singular) optimal dividend policy
  for the PDCP risk model. In this case, $D(t)$, the total amount of dividends paid up to time $t$,
is not necessarily absolutely continuous with respect to $t$.

Recall that for a given admissible dividend strategy $D=\{D(t),t\ge 0\}$, the associated EPV is given by \eqref{performance-defn} and the goal is to find an admissible
dividend strategy $D^*=\{D^*(t),t\ge 0\}$ that achieves the value function $V$ given by \eqref{value-defn}.
The following proposition indicates that $V$
 % defined in \eqref{value-defn}
is nondecreasing. It can be proved using exactly the same arguments as those used in \cite{Song-S-Z}.

\begin{prop}\label{prop1-sec4}
For any  $0\le y \le  x$, we have $V(x)  \ge   ( x-y) + V(y).
%\\ \label{value-increase-not-much}V(x) & \le V(y) +
%  V(x-y).
$
\end{prop}
Standard arguments using the dynamic programming principle and It\^o's formula lead to the following verification theorem, which enables us to identify the value function and an optimal dividend policy later. 
 % We leave the proof of Theorem \ref{thm-verification} to the supplement posted online. 

 \begin{thm}\label{thm-verification}
 Suppose there exists a  function $\vphi: \rr\mapsto \rr_+ \in C(\rr)\cap C^1(\rr \backslash{\mathbb D})$ with $  \vphi'(x+) < \infty$, $\vphi'(x-) < \infty$ for all $x\in {\mathbb D}$, where ${\mathbb D}$ is a countable set of points. Suppose $\vphi$    satisfies  $\vphi(y)=0$ for $y<0$ and that it  solves the following  quasi-variational inequality:
\beq\label{quasi-variational-ineq}
\max\set{(\op-\delta )\vphi(x), 1- \vphi'(x)}=0,\ \ x>0,
\eeq where $\op$ is the infinitesimal generator  defined  by   
% By virtue of \cite{Davis}, the generator of the PDCP is defined as
\beq\label{generator}
{\cal A} \vphi(x)= g(x) \vphi'(x) -\la \vphi(x)+ \la \int_0^\infty  \vphi(x-y)\df Q(y), \ \ x \ge 0,
\eeq
 % where $h$ is continuously differentiable.
\begin{itemize}\item[{\em (a)}] Then $\vphi(x) \ge V(x)$ for  every $x \ge 0$.

\item[{\em (b)}] Define the {\em continuation region}
 $ \mathcal C= \set{x\ge 0:
 1- \vphi'(x) < 0}. $
Assume
 there exists a dividend payment scheme $ \pi^*=\set{D^*(t): t\ge 0} \in \Pi$ and corresponding process $  X^*$ satisfying \eqref{controlled-surplus} such that,
 \begin{align}\label{cond-op0} &   X^*(t) \in \bar{\mathcal  C} \text { for Lebesgue almost all }  0 \le t \le \tau,   \\
 \label{cond-op1}
& \int_0^t \left[\vphi'(  X^*(s)) -1 \right] \df   D^{*c}(s)  =0, \text{ for any } t \le \tau, \\
\label{cond-op3}&  \lim_{N\to \infty} \ex_{x} \left[e^{-r(\tau\wedge N )} \vphi(  X^*(\tau\wedge N))\right]=0,
\end{align}
where  $D^{*c}(t):=D^*(t)- \sum_{0\le s \le t} \Delta D^{*c}(s)$ denotes the continuous part of $D^*$,
 and if $ X^*(s) \not=  X^*(s-)$, then
\begin{equation}\label{cond-op2}
  \vphi(  X^*(s))-\vphi( X^*(s-)) = -  \Delta D^*(s).
\end{equation}
Then $\vphi(x)= V(x)$ for every $x\ge 0$ and $\pi^* $ is an optimal dividend payment strategy.
\end{itemize}\end{thm}

% \begin{rem}\label{rem1-sec4}
% The conditions of \thmref{thm-verification} can be weakened. In fact, by virtue of
% \cite[Appendix D]{Oksendal}, we need only to assume that (i) $\vphi(\cdot,\al) \in C(\rr)\cap C^1(\rr \backslash{\mathbb D})$, 
% where ${\mathbb D}$ is a countable set of points, and (ii) $  \vphi'(x+) < \infty$, $\vphi'(x-) < \infty$ for all $x\in {\mathbb D}$.
% Under these conditions, there exist sequences $\set{\vphi_j}_{j=1}^\infty$, such that $\vphi_j\in
% C^1(\rr) $. Moreover, the following are satisfied:
% \begin{itemize}
% \item[{(a)}]  $\lim_{j\to \infty}\vphi_j(\cdot) \to \vphi(\cdot)$ uniformly on compact subsets of $\rr$,
% \item[(b)]  $\lim_{j\to \infty}(\op -\delta )\vphi_j(x) \to (\op -\delta)\vphi(x)$ 
% uniformly on compact subsets of $\rr\backslash{\mathbb D}$,  and
% \item[�] $\set{(\op-\delta)\vphi_j}_{j=1}^\infty$ is locally bounded on $\rr$.
% \end{itemize}
% Then, we can first work with the sequence $\phi_j$ in compact regions exactly the same
%  way as in the proof of \thmref{thm-verification}. Next, using (a), (b), and (c), 
% we can pass to the limit as $j\to \infty$ to obtain the same conclusions. 
% The reader is referred to   \cite{Oksendal} for details.
% \end{rem}

% Thanks to Remark \ref{rem1-sec4}, 
Using exactly the same arguments as those in the proof of Theorem \ref{thm-verification}, part (a),  we
obtain the following proposition.
\begin{prop}\label{prop2-sec4}
Suppose that there is function $\phi\in C(\rr)\cap C^1(\rr\backslash{\mathbb D})$, where ${\mathbb D}$ is a countable set,  satisfying  $\phi(y)=0$ for $y< 0$ and $\phi(x) \ge 0$, $\phi'(x) \ge \kappa >0$ for all $x\ge 0$. Then for any $x\ge 0$,
\beq\label{V-up-bd}V(x) \le \frac{1}{\kappa}\phi(x) + \frac{1}{\kappa}\sup_{D\in \Pi} \ex_{x} \int_0^\tau e^{-\delta s} (\op-\delta) \phi(X^D(s))\df s. \eeq 
\end{prop}
\begin{rem}
It follows that if there is a function $\phi$ satisfying the conditions of Proposition \ref{prop2-sec4} and also $\sup_{D\in \Pi} \ex_{x} \int_0^\tau e^{-rs} (\op-\delta) \phi(X^D(s))\df s < \infty$, then $V(x) < \infty$ for all $x\ge 0$. For example, if $g(x)=\rho x+ c$ for positive constants $\rho \le  \delta$ and $c$, then $V(x) < \infty$ for all $x\ge 0$. In fact, the function $\phi (x)=x^{+}, x\in \rr$,  satisfies the conditions of Proposition \ref{prop2-sec4}. Moreover, we compute for $x\ge 0$
\bed \begin{aligned}
(\op-\delta) \phi(x)  & = \rho x+ c - \lambda x +\lambda  \int_{0}^{\infty } (x-y) ^{+}\df Q(y) -\delta x \\
  &  \le (\rho-\delta)x +c -\lambda x + \lambda x  \int_{0}^{x} \df Q(y) - \lambda \int_{0}^{x} y \df Q(y)  \le c.
\end{aligned}\eed
 Then it follows that for any $x\ge 0$, $D\in \Pi$, we have $\ex_{x}\int_{0}^{\tau} e^{-\delta s} (\op-\delta) \phi(X(s))\df s \le \frac{c}{\delta} < \infty. $ Thus Proposition \ref{prop2-sec4} implies that  $V(x) < \infty$ for all $x\ge 0$.
\end{rem}

\begin{rem}\label{Rem-ViscositySoln}
In general the value function $V$ defined in \eqref{value-defn} is not necessarily smooth. Nevertheless, 
one can follow the arguments in \cite{Albrecher-Thon-08} to show that if $V$ is finite, then it is the unique {\em viscosity solution} of the quasi-variational inequality \eqref{quasi-variational-ineq}. 
\end{rem}

\subsection{Exponential Claims}
In order to obtain an explicit solution to the quasi-variational inequality
\eqref{quasi-variational-ineq} and an optimal dividend payment policy, as in Section \ref{sect-exponential1}, we again assume that the  claim sizes are exponentially distributed with mean $1/\al$ for some
  $\al>0 $.
% Then for any function $\vphi(x)$,
% \bed h(x):=\int_0^x \vphi(x-y) \df Q(y) = \int_0^x \vphi(x-y) \al e^{-\al y}
% \df y =\al e^{-\al x} \int_0^x \vphi(y) e^{\al y}\df y. \eed
% Hence it follows that
% \bed \(\frac{\df }{\df x}+ \al\)h(x) = \al \vphi (x). \eed
% Consequently, \eqref{quasi-variational-ineq}
% can be rewritten as
% \beq\label{qvi-differential-form}
% \max\set{g(x) \vphi''(x) + [\al g(x) + g'(x) -(\la +\delta)]
% \vphi'(x) -\al \delta \vphi(x), 1-\vphi'(x)} =0, \ \ x >0.
% \eeq
In such a case, % In what follows, we % will
we first construct an explicit solution to \eqref{quasi-variational-ineq}, which is exactly  the value function defined in \eqref{value-defn}.
 Then we provide easily verifiable sufficient conditions for the optimality of the barrier strategy.

%\nd{\bf Hypothesis B.}
%The integral-differential equation
%\beq\label{A-dl-equation}  (\op -\delta) \vphi(x) = g(x) \vphi'(x) -(\la+ \delta) \vphi(x) + \la \int_0^x \vphi(x-y) \al e^{-\al y}\df y =0, \quad x > 0,\eeq
% has a continuously differentiable and strictly increasing solution
% $\psi_1(x)$. Moreover,
% $\psi_1'(x)$ achieves its minimum value at $b > 0$ and
% $\psi_1'(x)$ is nondecreasing on $(b,\infty)$.

\begin{thm}\label{thm-soln-qvi} Suppose that $\psi_1$ is a positive and strictly increasing solution of the IDE \eqref{intdiff} and $\psi_1'(x)$ achieves its minimum value at $b > 0$ and
$\psi_1'(x)$ is nondecreasing on $(b,\infty)$. Then
\begin{itemize}
\item[\em (a)] the solution to \eqref{quasi-variational-ineq} is given by
\beq\label{soln-qvi} \Phi(x) =\begin{cases}\displaystyle \frac{\psi_1(x)}{\psi_1'(b)}, & \text{ if }\  0 \le x\le b,\\
\displaystyle x-b +  \frac{\psi_1(b)}{\psi_1'(b)}, & \text{ if }\ x> b.\end{cases}\eeq 
 \item[\em (b)] the {\em barrier strategy} given by continuous part 
 \begin{equation}
\label{barrier}
\df D^\ast(t)  =  g(b)I_{\set{X^*_t=b}} \df t ,
\end{equation}   
% \beq \label{barrier} \df D^\ast(t) =\begin{cases} 0, & \text{ if }\ X_t < b,\\
%  g(b)\df t, & \text{ if }\ X_t=b,\end{cases}\eeq 
and singular part \begin{equation}
\label{barrier-eq2}
\Delta D^\ast(t)= X^*_t-b, \qquad \mbox{ if }\ X^*_t > b, 
\end{equation} 
with $D^\ast(0-)=0$ is an optimal control that corresponds to $\Phi(x)$ given in \eqref{soln-qvi}, that is, $V(x)= \Phi (x) = J(x,D^*)$ for all $x\ge 0$.
\end{itemize}
\end{thm}

\begin{proof} (a)  Note that $\Phi \in C^1([0,\infty))$.
Obviously, if $x\le b$, $\Phi(x)$ satisfies \eqref{quasi-variational-ineq}.
If $x>b$,   $\Phi'(x) =1$. Therefore it remains to show that
\beq\label{ineq-claim0} g(x) \Phi'(x)  -(\la +\delta)\Phi (x) + \lambda \int_0^x \Phi(x-y) \al e^{-\al y}\df y \le 0, \qquad x > b. \eeq
To this end, we claim that
\beq\label{eq-claim1} g(x) \Phi''(x) + [\al g(x) + g'(x) -(\la +\delta)] \Phi'(x) -\al \delta \Phi(x) \le 0, \qquad x > b.\eeq
By assumption, $\psi_1'(x) >0$ for $ x>0$ and $\psi_1''(x) \ge 0$ for $x>b$. Hence it follows that for $x> b$
\beq\label{preparation-eq1} \begin{aligned}
g&(x) \Phi''(x) + [\al g(x) + g'(x) -(\la +\delta)] \Phi'(x) -\al \delta \Phi(x) \\
& \le g(x)\cdot \frac{\psi_1''(x)}{\psi_1'(x)} + [\al g(x) + g'(x) -(\la +\delta)] \frac{\psi_1'(x)}{\psi_1'(x)} -\al \delta \(x-b+ \frac{\psi_1 (b)}{\psi_1'(b)}\).
\end{aligned}\eeq
But $\psi_1'(x) $ is nondecreasing on $(b,\infty)$, thus we have
\beq\label{preparation-eq2} x-b = \int_b^x \frac{1}{\psi_1'(y)}\psi_1'(y) dy \ge \frac{1}{\psi_1'(x)}\int_b^x  \psi_1'(y) dy = \frac{1}{\psi_1'(x)}(\psi_1(x) -\psi_1(b)). \eeq
Since $\psi_1$ is a solution to \eqref{intdiff},
by applying the operator $(\df/\df x + \al)$ to \eqref{intdiff}, we see by straightforward calculations that
\beq\label{preparation-eq3}
g(x)\psi_1''(x) + [\al g(x) + g'(x) -(\la +\delta)] \psi_1'(x) -\al\dl \psi_1(x) =0. \eeq
A combination of \eqref{preparation-eq1}--\eqref{preparation-eq3} leads to
\bed \begin{aligned}
g&(x) \Phi''(x) + [\al g(x) + g'(x) -(\la +\delta)] \Phi'(x) -\al \delta \Phi(x) \\
 & \le \frac{1}{\psi_1'(x)}\left[g(x)\psi_1''(x) + [\al g(x) + g'(x) -(\la +\delta)] \psi_1'(x) -\al\dl \psi_1(x)\right] + \al\dl \psi_1 (b)\!\(\frac{1}{\psi_1'(x)}- \frac{1}{\psi_1'(b)}\)\\
 & =0 + \al\dl \psi_1 (b)\(\frac{1}{\psi_1'(x)}- \frac{1}{\psi_1'(b)}\)  \le 0,\end{aligned} \eed
using the fact
 % that $\psi_1 $ is solution to \eqref{A-dl-equation} and
that $\psi_1'(x) $ is nondecreasing on $(b,\infty)$.
Equation \eqref{eq-claim1} is therefore established.
 % and hence the function $\Phi$ satisfies \eqref{qvi-differential-form}.

Next we show that $\Phi$ satisfies \eqref{ineq-claim0}.
It follows immediately from \eqref{eq-claim1} that 
\beq \label{diffclaim}\begin{aligned}
&\frac{\df}{\df x} \left[e^{\alpha x} \left( g(x)\Phi'(x)-(\lambda+\delta)\Phi(x)+\lambda \int^x_0 \Phi(x-y) \alpha e^{-\alpha y} \df y  \right)  \right]\\
&=e^{\alpha x} \Big(g(x) \Phi''(x) + [\al g(x) + g'(x) -(\la +\delta)] \Phi'(x) -\al \delta \Phi(x)\Big) \le 0.
\end{aligned}
\eeq
Note that $ g(b) \Phi'(b) - (\la+\delta) \Phi(b) + \la \int^b_0 \Phi(b-y) \alpha e^{-\alpha y} \df y=0.$
Integrating both sides of the inequality \eqref{diffclaim} proves the claim \eqref{ineq-claim0}.

(b) It is easy to verify that the strategy $D^*$ and the corresponding
surplus process $X^*$ satisfy all the conditions in \thmref{thm-verification}(b).
Hence part (a) of this theorem and Theorems \ref{thm-verification}   imply that $J(x,D^*)= \Psi(x)= V(x)$ for all $x\ge 0$.
\end{proof}

	As in Section \ref{sect-exponential1},  the next theorem provides easily verifiable conditions under which the optimality of the barrier strategy is established. Our result is applicable in a general PDCP risk model with minimal assumptions.
\begin{thm} \label{finalbar}
 % Suppose that
% $\displaystyle \sup_{x \ge 0} \left\{g''(x)+\alpha g'(x)-\alpha \delta \right\}<0.$
Assume \eqref{eq:sup-g-x}.
Then the following assertions are valid. \begin{enumerate}
\item[(i)] If  %$\alpha \lambda g(0)+(g'(0)-\lambda-\delta)(\lambda+\delta)>0$,
in addition, \eqref{eq:g-0}  is also satisfied,
    then the  unrestricted payment scheme given by   \eqref{barrier}--\eqref{barrier-eq2} with barrier level $b>0$ is optimal and the value function is \eqref{soln-qvi}.
\item[(ii)] Otherwise, if \eqref{eq:g-0}  is not satisfied, that is,  \begin{equation}
\label{eq:g<=0}
\alpha \lambda g(0)+(g'(0)-\lambda-\delta)(\lambda+\delta)\le 0,
\end{equation}
 then the optimal unrestricted payment scheme is % defined by $\Delta D^\ast(0)=x$ and $\df D^\ast (t)=g(0) \df t$ if $X_t=0$ 
 given by \eqref{barrier}--\eqref{barrier-eq2} with the barrier level $b=0$ and the value function is given by
\[V(x)=x+\frac{g(0)}{\lambda+\delta},\qquad x \ge 0.\]
\end{enumerate}
\end{thm}
\begin{proof}
(i)
We have shown in the proof of Theorem \ref{finalthr} that \eqref{intdiff} admits a positive and increasing solution $\psi_1$ satisfying $\psi_1''(b)=0$, $\psi_1''(x) < 0$ for $x\in [0,b)$ and $\psi_1''(x) > 0$ for $x\in (b,\infty)$, where $b >0$.
Therefore $\psi_1$ satisfies the conditions in Theorem \ref{thm-soln-qvi} and hence the desired assertion follows.
% from Theorems \ref{thm-soln-qvi} and
% \ref{thm-barrier}.

(ii) In this case, \eqref{eq:g<=0} implies that $\psi''_1(0) \ge 0$.
% As shown in the first part of the proof of Lemma \ref{idecond},
Note that $\psi'_1(x)>0$ for all $x \ge 0$.
Differentiating \eqref{eq-psiode} we obtain
\[g(x)\psi'''_1(x)+(2g'(x)+\alpha g(x)-\lambda-\delta)\psi''_1(x)+(g''(x)+\alpha g'(x)-\alpha \delta) \psi'_1(x)=0.\] Since $g''(x)+\alpha g'(x)-\alpha \delta <0$ for all $x \ge 0$, we can again show by contradiction that $\psi'_1(x)$ cannot have a local maximum, which implies that $\psi''_1(x) >0$ for all $x> 0$. Following the proof of \eqref{ineq-claim0} with $\Phi(x)=x+\psi_1(0)/\psi'_1(0)=x+g(0)/(\lambda+\delta)$ and $b=0$, we can show that $(\mathcal{A}-\delta)\Phi(x) \le 0$ for all $x >0$. Thus $\Phi(x)$ clearly solves the quasi-variational inequality \eqref{quasi-variational-ineq} and consequently $\Phi(x) \ge V_R(x)$ for every $x \ge 0$. 
Last, we can %follow the proof of \eqref{x>b} with $b=0$ to show 
easily verify that the stated dividend scheme achieves the upper bound $\Phi(x)$ and hence is the optimal unrestricted policy.
\end{proof}

\begin{rem}
In the complete generality of the PDCP model, we have shown that under the assumption \eqref{eq:sup-g-x},
 % $\sup_{x \ge 0} \left\{g''(x)+\alpha g'(x)-\alpha \delta \right\}<0$
  the optimal unrestricted dividend scheme is always the barrier strategy, which is to pay out the initial surplus as dividends in excess of a certain barrier level and then pay dividends continuously at the rate of all incoming cash flow so as to keep the surplus at the barrier until the time of ruin. Under the condition
   % $\alpha \lambda g(0)+(g'(0)-\lambda-\delta)(\lambda+\delta)>0,$
   \eqref{eq:g-0}, the barrier level is chosen to be a positive level. Otherwise, if \eqref{eq:g<=0} holds, then the barrier level is set to be zero.
\end{rem}

\section{Examples} \label{examples}

 The following examples demonstrate the general results developed in Sections  \ref{sect-restricted} and \ref{sect-singular}.
Compared with  the existing literature (\cite{Alb-Har-2007,Fang2007,Gerber2006}, and \cite{schmidli2008}), our solution method to optimal dividend policies
 is simpler yet has an advantage of being applicable to general PDCP risk models. The last example with regressive growth rate appears to be new in the literature.

 %\subsection{General Case} \label{general}

\subsection{Cram\'er-Lundberg Model} \label{clmodel}

To illustrate our results, let us consider the special case when $g(x)\equiv c>0$ and the claim size distribution
$ Q(y)=1-e^{-\alpha y}$,  $ y\geq 0. $  Note that the  condition  \eqref{eq:g-0} in Theorem \ref{finalthr} is in fact
\begin{equation}\label{eq-b>0}\al\la c> (\la+\delta)^2,\end{equation}
and the second condition \eqref{eq:sup-g-x} is automatically satisfied.
One simply has to find a solution $\psi_1(x)$ to the integro-differential equation
\beq c\varphi'(x)-(\lambda+\delta)\varphi(x)+\lambda \int^x_0 \varphi(x-y) \alpha e^{-\alpha y} \df y=0, \qquad x>0,\label{expintdiff}\eeq and an increasing and concave solution $\psi_2(x)$ to the  differential equation
\beq (c-u_0)\varphi''(x)-[\alpha c-\alpha u_0-(\lambda+\delta)]\vphi'(x)-\alpha \delta \vphi(x)=0, \qquad x>0.\label{expintdife2}\eeq
The unique solution (up to a constant multiple) to
\eqref{expintdiff} is
$ \psi_1 (x)= (r+\al)e^{rx}-(s+\al)e^{sx},  $ 
where $-\al<s<0< r$ are the roots of
\beq\label{eq-char-eq} c \xi^2 -(\la+\delta -\al c)\xi -\al\delta=0.\eeq Note that $\psi_1$ is positive and strictly increasing on $[0,\infty)$.
Similarly, the differential equation {\eqref{expintdife2}} has an increasing and
%strictly increasing
 concave solution
$ \psi_2(x)= - e^{tx}, $ where $t$ is the negative root of
\beq\label{eq-char-eq1} (c-u_0) \xi^2 -(\la+\delta -\al c+\al u_0)\xi -\al\delta=0.\eeq

\subsubsection{Restricted Payment Scheme}
By  virtue of condition \eqref{condn}, we obtain
\[\frac{(r+\alpha)e^{rd}-(s+\alpha) e^{sd}}{r(r+\alpha)e^{rd}-s(s+\alpha) e^{sd}}=\frac{1}{t}+\frac{u_0}{\delta}.\]
 %It follows immediately that
Solving the above equation for $d$ gives
\beq\label{d-threshold} d=\frac{1}{r-s} \ln \left[ \frac{(s+\alpha)(\delta t-\delta s- st u_0)}{(r+\alpha)(\delta t-\delta r- rt u_0)} \right]= \frac{1}{r-s} \ln \left[ \frac{s(s-t)}{r(r-t)}\right] ,\eeq
which agrees with (9.15) of \cite{Gerber2006}. %(2006).
However, our approach is considerably simpler than their  method of optimizations.
%used in Gerber and Shiu (2006).
Note that the following are all equivalent sufficient and necessary condition for $d>0$
\beq \frac{1}{t}> \frac{c}{\lambda+\delta}-\frac{u_0}{\delta}, \mbox{ or }\, \frac{1}{t}+\frac{u_0}{\delta} >\frac{r-s}{r(r+\alpha)-s(s+\alpha)}.\label{b0}\eeq

Therefore, according to Theorem \ref{finalthr}(i), if \eqref{b0} is satisfied,
then   the
value function is
\begin{equation*} V_R(x) =\begin{cases} \displaystyle \frac{(r+\al)e^{rx}-(s+\al)e^{sx}}{r(r+\al)e^{rd}-s(s+\al)e^{sd}}, & \text{ if }\ 0 \le x < d;\\
 \displaystyle \frac{u_0}{\delta} +  \frac{1}{t} e^{t(x-d)}, & \text{ if }\ x\ge d,\end{cases} \end{equation*}
and the optimal restricted dividend payment scheme is the threshold strategy given in \eqref{optimal-restricted-policy}.
Otherwise, if \eqref{b0} is not satisfied, then according to Theorem \ref{finalthr}(ii), the value function is
\[V_R(x)=\frac{u_0}{\delta}\left[1-\frac{\lambda}{\lambda+\delta-(c-u_0)t} e^{-tx}  \right],\qquad x \ge 0,\] and the optimal restricted dividend policy is to pay at the maximal rate $u_0$ at all time until ruin.

%Assume that $d>0$, or equivalently,
%\beq\label{exm1-sec3-cond1}  \frac{(s+\alpha)(\delta t-\delta s- st u_0)}{(r+\alpha)(\delta t-\delta r- rt u_0)} >1.\eeq
%We claim that  $\psi_1 $   is concave on the interval  $(0,d)$.
%   In fact, it is straightforward to verify that the function $\psi_1'(x)=r(r+\al)e^{rx}-s(s+\al)e^{sx}$ is decreasing on $(-\infty, b)$ and increasing on $[b,\infty)$, where
% \beq\label{b-barrier} b= \frac{1}{r-s}\ln \(\frac{s^2(s+\al)}{r^2(r+\al)}\) =\frac{1}{r-s}\ln \(\frac{s[(\la+\delta)s + \al\delta]}{r[(\la+\delta)r + \al\delta]}\).
%\eeq
%Therefore the desired concavity will follow if we can show that $d \le b$.
% Recall that $-\al < s < 0 < r$. Thus a comparison between
%\eqref{d-threshold} and \eqref{b-barrier} reveals that it suffices to prove
%\beq\label{eq-d<b} \frac{s(s-t)}{r(r-t)}- \frac{s[(\la+\delta)s + \al\delta]}{r[(\la+\delta)r + \al\delta]} < 0,\eeq
%or, equivalently,
%\beq\label{eq-d<b-2} \frac{(s-t)}{(r-t)}- \frac{[(\la+\delta)s + \al\delta]}{ [(\la+\delta)r + \al\delta]}=\frac{(s-r)[\al\delta+ t(\lambda+ \delta)]}{(r-t)[(\la+\delta)r + \al\delta]} > 0.\eeq
%But since $$t=\frac{\la+\delta-\al(c-u_0) -\sqrt{(\la+\delta-\al(c-u_0))^2+ 4(c-u_0)\al\delta}}{2(c-u_0)},$$
%straightforward algebra reveals that $\al\delta+ t(\lambda+ \delta) < 0$, which, in turn, leads to \eqref{eq-d<b-2}. Therefore it follows that $\psi_1'(x)$ is decreasing on $(0,d)$ and hence $\psi_1(x)$ is indeed concave on $(0,d)$.

\subsubsection{Unrestricted Payment Scheme}
% Recall that $\psi_1(x)=(r+\al)e^{rx}- (s+\al) e^{sx} $ solves the integral-differential equation \eqref{expintdiff} and
Note that $\psi'_1(x)=r(r+\al)e^{rx}- s(s+\al) e^{sx} $ achieves its unique minimum value at
\bed b= \frac{1}{r-s}\ln \(\frac{s^2(s+\al)}{r^2(r+\al)}\) =\frac{1}{r-s}\ln \(\frac{s[(\la+\delta)s + \al\delta]}{r[(\la+\delta)r + \al\delta]}\),
\eed
and that $\psi'$ is nondecreasing on $(b,\infty)$.
Therefore in view of Theorem \ref{finalbar}(i), under the condition \eqref{eq-b>0}, $b>0$ and the dividend payment strategy defined in \eqref{barrier} is optimal and the
value function is
\bed V(x)= \begin{cases}\displaystyle \frac{(r+\al)e^{rx}-(s+\al)e^{sx}}{r(r+\al)e^{rb}- s(s+\al) e^{sb}} &  \quad  \text{ if } x< b,\\[.5ex]
\displaystyle  x-b+  \frac{(r+\al)e^{rb}-(s+\al)e^{sb}}{r(r+\al)e^{rb}- s(s+\al) e^{sb}}  & \quad \text{ if } x\ge  b.\end{cases}\eed
Although unnecessary as the result is proved in more generality, we can verify using the fact that $r$ and $s$ are the roots of \eqref{eq-char-eq} that $ b> 0$  if and only if \eqref{eq-b>0} holds.

 On the other hand, according to Theorem \ref{finalbar}(ii), % under the condition
%  \[\al\la c \le (\la+\delta)^2, \]
if \eqref{eq-b>0} is not satisfied,  the barrier level is set to zero. Hence it follows that
  $$V(x)= x+ \frac{c}{\la+\delta},\qquad x \ge 0.$$

 Therefore we can summarize the value function as
 \bed V(x)= \begin{cases}\displaystyle \frac{(r+\al)e^{rx}-(s+\al)e^{sx}}{r(r+\al)e^{rb}- s(s+\al) e^{sb}} &  \quad  \text{ if } \al\la c> (\la+\delta)^2 \text{ and } x< b,\\[.5ex]
\displaystyle x-b+  \frac{(r+\al)e^{rb}-(s+\al)e^{sb}}{r(r+\al)e^{rb}- s(s+\al) e^{sb}}  & \quad \text{ if } \al\la c> (\la+\delta)^2 \text{ and } x\ge  b, \\[.5ex]
 \displaystyle x+ \frac{c}{\la+\delta} & \quad  \text{ if } \al\la c\le  (\la+\delta)^2,\end{cases}\eed
which agrees results from \cite[p. 94]{schmidli2008}. But our approach
 is much simpler than theirs.

\begin{figure}[htbp]
\begin{center}
 \mbox{\subfigure[Value functions $V(x)$ and $V_R(x)$]{\epsfig{figure=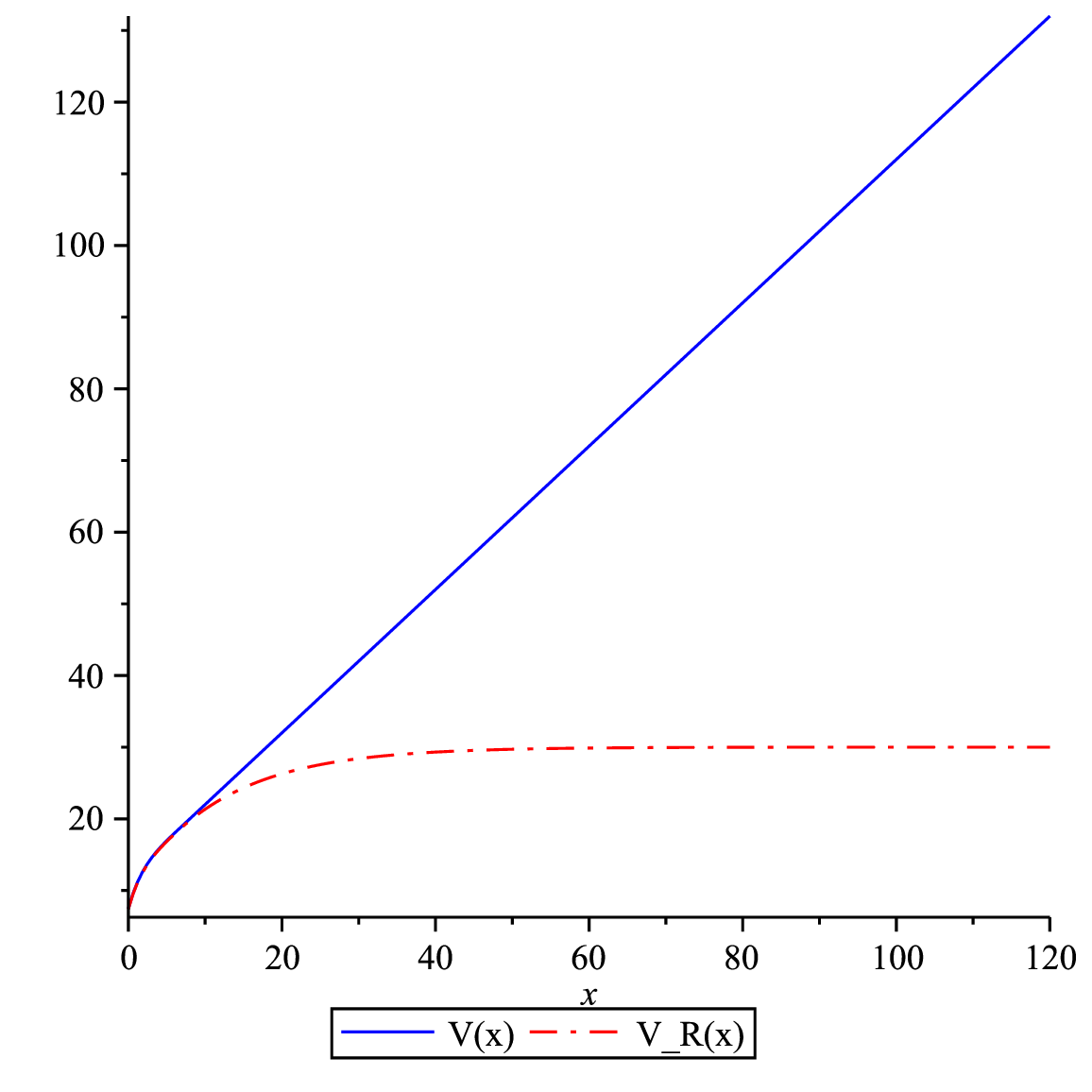,width=0.43\linewidth}}}\ \
 \mbox{\subfigure[The difference $V(x)-V_R(x)$]{\epsfig{figure=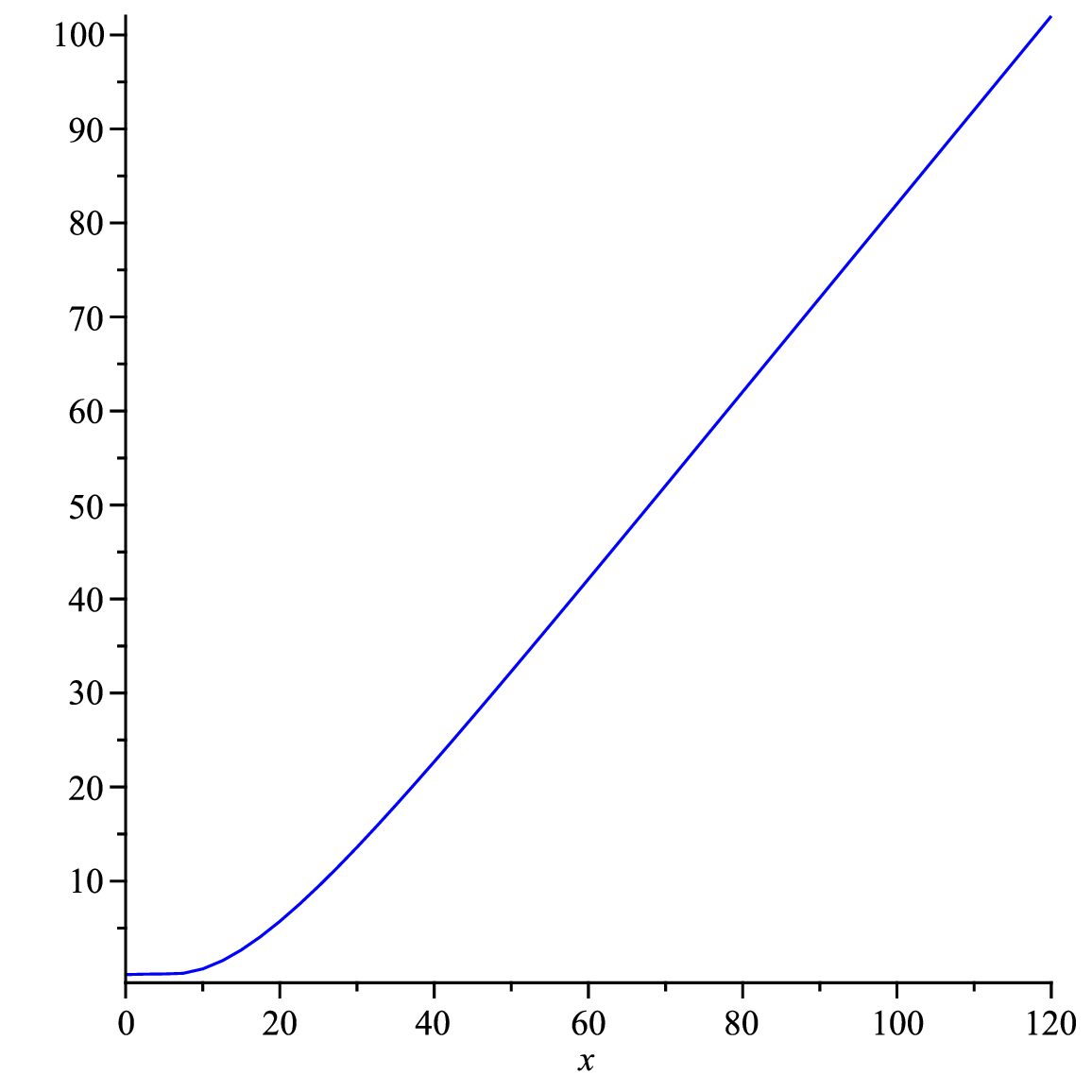,width=0.43\linewidth}}}
\caption{Comparison %of optimal policies
in Cram\'er-Lundberg model} \label{fig1}
\end{center}
\end{figure}

Finally we demonstrate the comparison of restricted and unrestricted payment schemes through a numerical example, in which $\al=1$, $\delta=0.1$, $c=4$, $\lambda=2$, and $u_0=3$. Note that both \eqref{eq-b>0} and \eqref{b0} are satisfied and the threshold and barrier levels are determined to be $ d=6.291707724$ and  $b=7.004651047$, respectively. 
The resulting unrestricted and restricted value functions $V(x)$ and $V_R(x)$ are shown in Figure \ref{fig1}(a).
We also plot the difference $V(x)-V_R(x)$ in Figure \ref{fig1}(b).
Note that  %$V(x)>V_R(x)$ for all $x\ge 0$ in this example. Also,
 the plot of $V_R(x)$ in Figure \ref{fig1}(a)
  %Finally, to
also demonstrates the limit result of $V_R(x)$ stated in Lemma \ref{lem-bounds}.
  %we plot  $ V_R(x)-\frac{u_0}{\delta}$ in Figure \ref{fig2} (b).

\subsection{Classical Model with Constant Interest} \label{interestmodel}
In this case, we consider $g(x)=\rho x+c$, where $\rho,\,c$ are positive constants. Note that the first condition  \eqref{eq:g-0} in Theorem \ref{finalthr} is
\beq \label{intcond1}\alpha \lambda c+(\rho-\lambda-\delta)(\lambda+\delta)>0.\eeq And the second condition \eqref{eq:sup-g-x} is satisfied if and only if
\beq \rho <\delta .\label{intcond2}\eeq

We now need to find a solution $\psi_1(x)$ to the IDE
\beq \label{interestide}
(\rho x+c) \psi'(x)-(\lambda+\delta) \psi(x)+\lambda \int^x_0 \psi(x-y) \alpha e^{-\alpha y} \df y=0, \qquad x>0,
\eeq and an increasing and concave solution $\psi_2(x)$ to the ODE
\beq \label{interestode}
(\rho x+c-u_0) \psi''(x)+[\alpha (\rho x+c)-\alpha u_0+\rho-(\lambda+\delta)]\psi'(x)-\alpha \delta \psi(x)=0,\qquad x >0.
\eeq
Applying $(\df/\df x+\alpha)$ to both sides of \eqref{interestide} gives
\beq (\rho x+c)\psi''(x)+[\alpha(\rho x+c)+\rho-\lambda-\delta] \psi'(x)-\alpha \delta \psi(x)=0. \label{interestode1} \eeq
Then a fundamental system of \eqref{interestode1} is given by
\[p_1(x)=z^{b-1} e^{-z} M(a,b,z),\qquad p_2(x)=z^{b-1} e^{-z} U(a,b,z),\] where $M$ and $U$ are the Kummer's functions of the first and second kind respectively, and
\[z=\frac{\alpha}{\rho}(\rho x+c),\qquad a=1+\frac{\delta}{\rho}>0 ,\qquad b=1+\frac{\lambda+\delta}{\rho}>0.\]
Using \cite[p.325-326, (13.3.20),(13.3.27)]{Olver:2010:NHMF},
\bed \begin{aligned}
\frac{\df}{\df z}(e^{-z}M(a,b,z))&=(-1)\frac{b-a}{b}e^{-z}M(a,b+1,z),\quad
\frac{\df}{\df z}(e^{-z}U(a,b,z))&=(-1)e^{-z}U(a,b+1,z), \end{aligned}
\eed it is easy to show that
\bed \begin{aligned} cp'_1(0)-(\lambda+\delta)p_1(0)&=-\left(\frac{\alpha c}{\rho}\right)^b \frac{\lambda \rho}{\rho+\lambda+\delta} e^{-\alpha c/\rho}M\left(a,b+1,\frac{\alpha c}{\rho} \right),\\
cp'_2(0)-(\lambda+\delta)p_2(0)&=-\left(\frac{\alpha c}{\rho}\right)^b \rho e^{-\alpha c/\rho}U\left(a,b+1,\frac{\alpha c}{\rho} \right). \end{aligned}
\eed It follows from \eqref{interestide} that $c\psi'(0)-(\lambda+\delta)\psi(0)=0$. Hence the unique solution to \eqref{interestide} up to a multiplicative constant is given by
\[\psi_1(x)=K_1 p_1(x)- K_2 p_2(x).\]
where
\bed
K_1=(\rho+\lambda+\delta) U\left(a,b+1, \frac{\alpha c}{\rho}\right),\qquad K_2=\lambda M\left(a,b+1,\frac{\alpha c}{\rho} \right).
\eed
By \cite[p.326, (13.3.28)]{Olver:2010:NHMF},
\bed
%\begin{aligned}
%\frac{\df^n}{\df z^n}(z^{b-1}e^{-z}M(a,b,z))&=(b-n)_n z^{b-n-1}e^{-z}M(a-n,b-n,z);\\
\frac{\df^n}{\df z^n}(z^{b-1}e^{-z}U(a,b,z))=(-1)^n z^{b-n-1}e^{-z}U(a-n,b-n,z),
%\end{aligned}
\eed we see that an increasing and concave solution to \eqref{interestode} is given by
\[\psi_2(x)=-p^\ast_2(x)=-(z^\ast)^{b-1} e^{-z^\ast} U(a,b,z^\ast),\qquad \mbox{ where } z^\ast=\frac{\alpha}{\rho}(\rho x+c-u_0).\]

\subsubsection{Restricted Payment Scheme}

We want to determine the unique number $d>0$ for which equation \eqref{condn} holds. Among other equivalent conditions, one sufficient and necessary condition for $d>0$ is
\beq \displaystyle \frac{\alpha}{\rho}(c-u_0) \frac{U\left(1+\delta/\rho,1+(\lambda+\delta)/\rho,\alpha(c-u_0)/\rho  \right)}{U\left(\delta/\rho,(\lambda+\delta)/\rho,\alpha(c-u_0)/\rho  \right)} < \frac{u_0(\lambda+\delta)-c \delta}{(\lambda+\delta) \delta}.\label{intcond3}\eeq According to Theorem \ref{finalthr}(i), under the conditions \eqref{intcond1}, \eqref{intcond2}, \eqref{intcond3}, the value function is given by
\bed V_R(x)= \begin{cases}\displaystyle \frac{K_1\, p_1(x)-K_2\,p_2(x)}{K_1 C_1(d)+K_2 C_2(d)} &  \quad  \text{ if } 0\le x< d;\\[.5ex]
 \displaystyle \frac{u_0}{\delta} -\frac{p^\ast_2(x)}{C_3(d)}   & \quad \text{ if } x\ge  d,\end{cases}\eed
where $C_1(\cdot), C_2(\cdot), C_3(\cdot)$ are positive functions given by
\bed
\begin{aligned}
C_1(d) &=\frac{\lambda+\delta}{\rho} \left( \alpha d+ \frac{\alpha c}{\rho} \right)^{(\lambda+\delta)/\rho-1}\exp \left\{ -\alpha d+\frac{\alpha c}{\rho}\right\}M\left(\frac{\delta}{\rho},\frac{\lambda+\delta}{\rho},\alpha d+\frac{\alpha c}{\rho}\right);\\
C_2(d) &= \left( \alpha d+ \frac{\alpha c}{\rho} \right)^{(\lambda+\delta)/\rho-1}\exp \left\{ -\alpha d+\frac{\alpha c}{\rho}\right\}U\left(\frac{\delta}{\rho},\frac{\lambda+\delta}{\rho},\alpha d+\frac{\alpha c}{\rho}\right);\\
C_3(d) &=\left[ \alpha d+\frac{\alpha}{\rho}(c-u_0) \right]^{(\lambda+\delta)/\rho-1}\!\! \exp \left\{ -\alpha d +\frac{\alpha}{\rho}(c-u_0) \right\} U\left(\frac{\delta}{\rho},\frac{\lambda+\delta}{\rho}, \alpha d+\frac{\alpha}{\rho}(c-u_0) \right).
\end{aligned}
\eed
If \eqref{intcond3} is violated, it follows from Theorem \ref{finalthr}(ii) that the value function is given by
\[V_R(x)=\frac{u_0}{\delta}-\frac{p^\ast_2(x)}{K},\qquad x \ge 0,\]
where
%\bed \begin{aligned} K&=\left(\frac{\alpha(c-u_0)}{\rho}\right)^{(\lambda+\delta)/\rho}e^{-\alpha(c-u_0)/\rho}\\
%&\times \frac{(\rho/\alpha)U(\delta/\rho,(\lambda+\delta)/\rho,\alpha(c-u_0)/\rho)+(\lambda+\delta)U(\delta/\rho+1,(\lambda+\delta)/\rho+1,\alpha(c-u_0)/\rho}{(\lambda/\delta)u_0}\end{aligned}
%\eed
\bed \begin{aligned} K=& \Big(\frac{\alpha(c-u_0)}{\rho}\Big)^{1+(\lambda+\delta)/\rho}e^{-\alpha(c-u_0)/\rho}\left(\frac{\delta \rho}{\lambda u_0}\right)U\left(\frac{\delta}{\rho}+1,\frac{\lambda+\delta}{\rho}+2,\frac{\alpha(c-u_0)}{\rho}\right).\end{aligned}
\eed

\subsubsection{Unrestricted Payment Scheme}
Note that $\psi_1(x)$ solves the IDE \eqref{interestide}. Using \cite[p.325-326, (13.3.21), (13.3.28)]{Olver:2010:NHMF}, we can determine the unique minimum value of $\psi'_1(x)$ at $b$ by solving the equation
\bea  K_1\left(\frac{\lambda+\delta}{\rho}-1 \right) \ad \left( \frac{\lambda+\delta}{\rho}\right)M\left( \frac{\delta}{\rho}-1,\frac{\lambda+\delta}{\rho}-1,\alpha b+\frac{\alpha c}{\rho} \right)\\ \aad +K_2 U\left(\frac{\delta}{\rho}-1,\frac{\lambda+\delta}{\rho}-1,\alpha b+\frac{\alpha c}{\rho}\right)=0.\eea
According to Theorem \ref{finalbar}(i), under % the
conditions \eqref{intcond1} and \eqref{intcond2}, we must have $b>0$ and that the dividend payment policy defined in \eqref{barrier} is optimal and the value function is given by
\bed V(x)= \begin{cases}\displaystyle \frac{K_1\, p_1(x)-K_2 \, p_2(x)}{K_1 C_1(b)+K_2 C_2(b)} &  \quad  \text{ if } 0\le x< b;\\[.5ex]
 \displaystyle x-b + \frac{K_1\, p_1(b)-K_2 \, p_2(b)}{K_1 C_1(b)+K_2 C_2(b)}, & \quad \text{ if } x\ge  b.\end{cases}\eed
In the case when 
$ \alpha \lambda c+(\rho-\lambda-\delta)(\lambda+\delta)\le 0, $ 
Theorem \ref{finalbar}(ii) tells us that the value function of the optimal payment scheme is
\[V(x)=x+\frac{c}{\lambda+\delta},\qquad x \ge 0.\]

\begin{figure}[htbp]
\begin{center}
 \mbox{\subfigure[Value functions $V(x)$ and $V_R(x)$]{\epsfig{figure=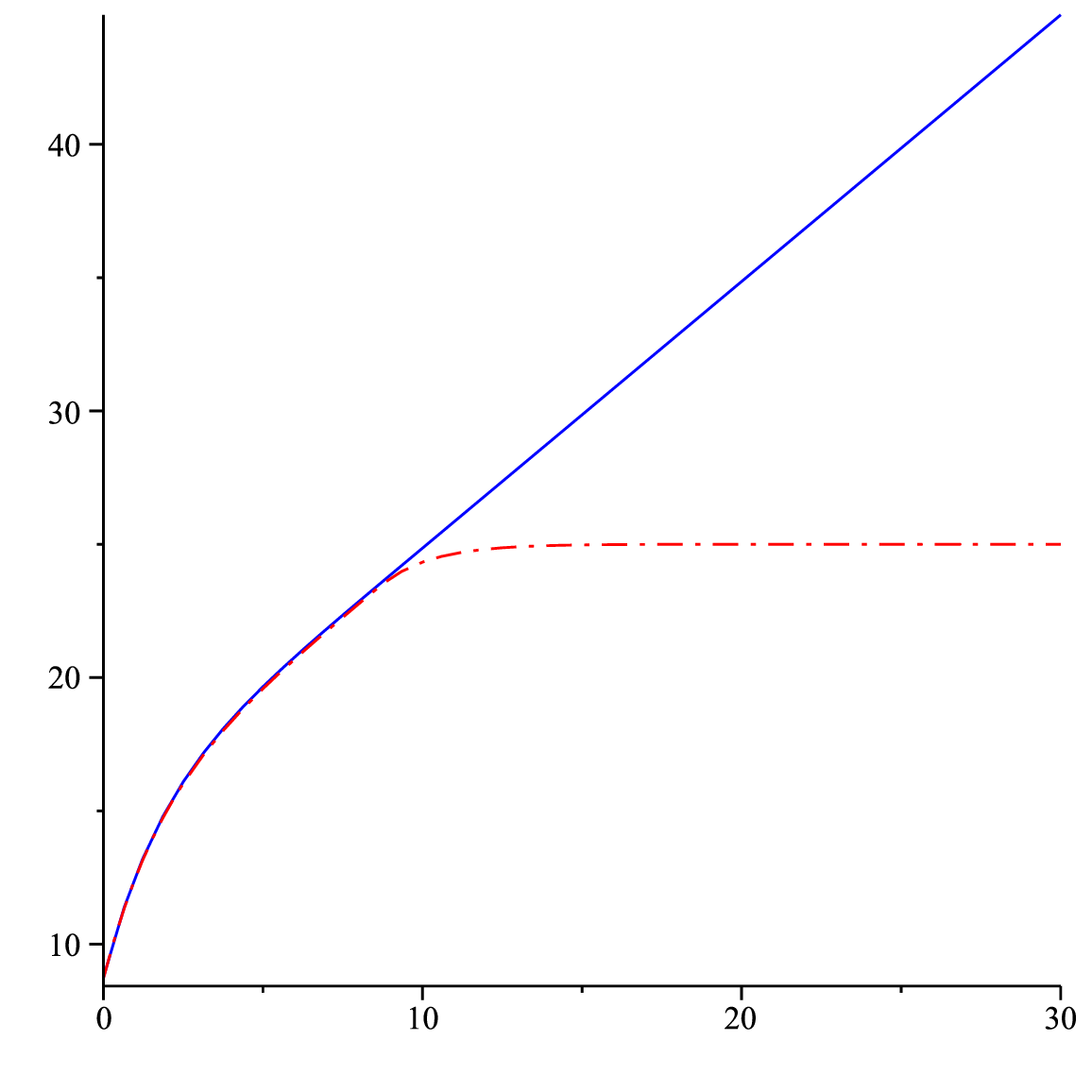,width=0.43\linewidth}}}\ \
 \mbox{\subfigure[The difference $V(x)-V_R(x)$]{\epsfig{figure=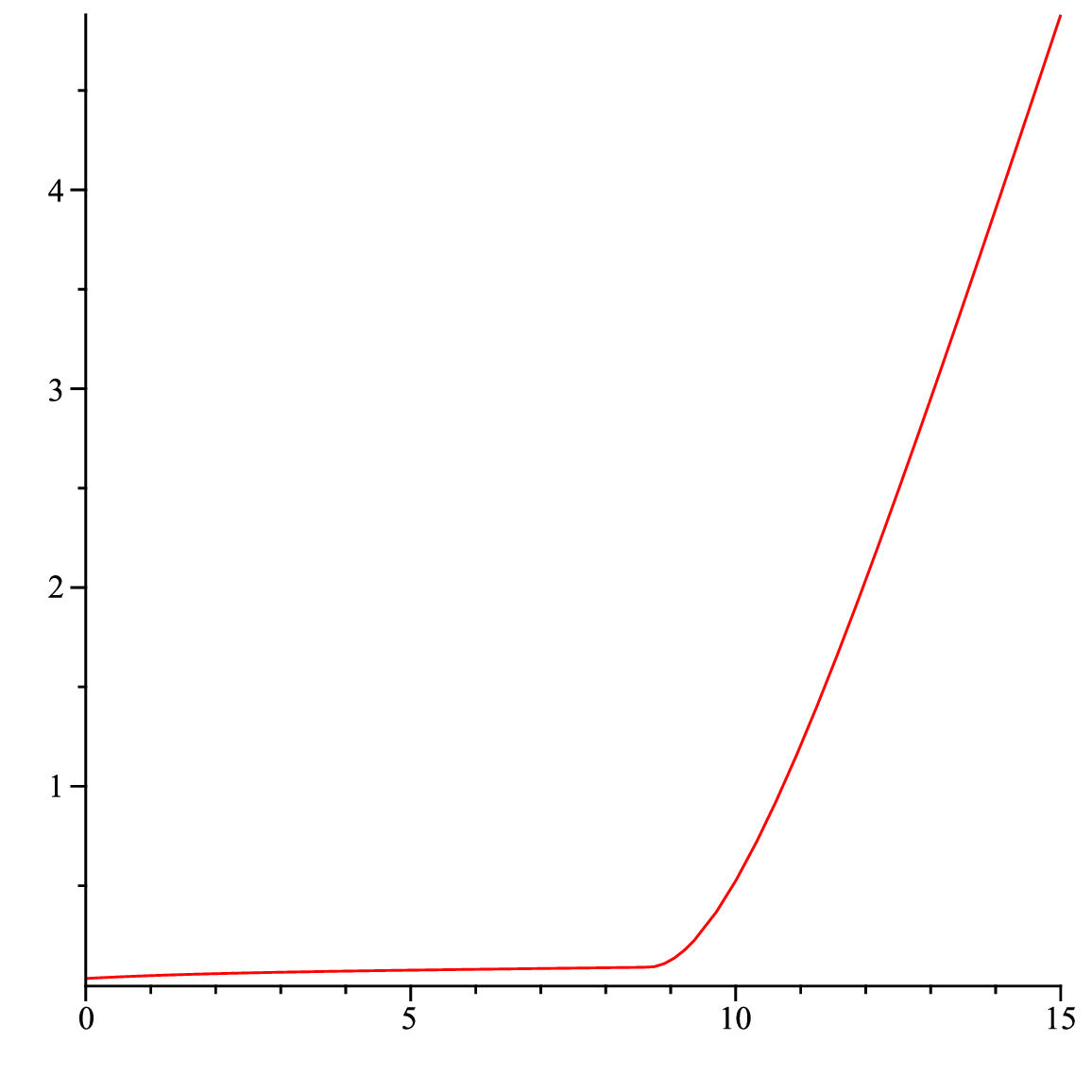,width=0.43\linewidth}}}
\caption{Comparison %of optimal policies
 in classical model with interest} \label{fig2}
\end{center}
\end{figure}

Last, we make a comparison of the optimal restricted and unrestricted payment schemes developed in this model with a numerical example. The parameters are chosen as follows: $\rho=0.05, \lambda=1, \delta=0.06, \alpha=1, c=2$, and $ u_0=1.5.$ Note that all three conditions \eqref{intcond1}, \eqref{intcond2} and \eqref{intcond3} are satisfied, which guarantees the existence and uniqueness of the constants $0<d<b$.
In this example, we calculate $d=8.660487436$ and $b=9.870353657.$ The comparison of the value functions of both optimal unrestricted and restricted dividend policies are plotted in Figure \ref{fig2}(a) and their difference shown in Figure \ref{fig2}(b).

\subsection{Regressive growth model}

We consider a risk model with a regressive growth rate. Assume that the growth rate $g$ of the overall insurance business depends on the insurance surplus and can be adjusted by changing the premium rates according to the following stochastic differential equation.  For $b>0$,
\[\df g(X_t)=b[c-g(X_t)] \df X_t.\]  When the growth rate is below a normal rate $c>0$, the incremental growth is accelerated. When the growth rate is above the normal rate $c$, the incremental growth is reduced. An application of   Ito's formula shows that
\[g(X_t)=a e^{-bX_t}+c,\qquad \mbox{ where } a=g(0)-c.\] Since a high premium is charged when the surplus is low,  we assume that $a>0$.

We are now interested in finding the optimal dividend payment schemes in the risk model \eqref{controlled-surplus} with $g(x)=a e^{-bx}+c$ where $a,b>0$ and $c> u_0\ge 0$. In this case, the first condition \eqref{eq:g-0} in Theorem \ref{finalthr} translates to
\be \alpha \lambda(a+c)-(ab+\lambda+\delta)(\lambda+\delta)>0. \label{newcond1}\ee The second condition \eqref{eq:sup-g-x} is true if and only if
\be ab (b-\alpha)<\alpha \delta. \label{newcond2}\ee
We look for a non-trivial solution $\psi_1$ to the IDE (which is positive and increasing if $\psi_1(0)>0$, and negative and decreasing if $\psi_1(0)<0$ by Lemma \ref{sec34-lem1})
\be
(ae^{-bx}+c)\psi'(x)-(\lambda+\delta)\psi(x)+\lambda \int^x_0 \psi(x-y) \alpha e^{-\alpha y} \df y=0,\qquad x >0, \label{newide}
\ee and an increasing and concave solution $\psi_2$ to the ODE
\be (ae^{-bx}+c-u_0) \psi''(x)+[a(\alpha-b)e^{-bx}+\alpha (c-u_0)-\lambda-\delta]\psi'(x)
-\alpha \delta \psi(x)=0. \label{newode}\ee Applying $(\df /\df x+\alpha)$ to \eqref{newide} and letting $z=(-c/a)e^{bx}, y(z)=\psi(x)$,  we arrive at the Gauss' hypergeometric equation
\[z(z-1) y''(z)+\left[(k_1+k_2+1)z-\frac{\alpha}b\right]y'(z)+k_1 k_2 y(z)=0,\] where 
\[k_1=\frac{\alpha c-\lambda-\delta+\sqrt{(\alpha c-\lambda-\delta)^2+4\alpha\delta c}}{2bc}, \qquad k_2=\frac{\alpha c-\lambda-\delta-\sqrt{(\alpha c-\lambda-\delta)^2+4\alpha\delta c}}{2bc}.\] The equation has regular singularities at $z=0,1,\infty$ and \cite[p.395]{Olver:2010:NHMF} give six solutions, among which the following two will be used:
\[F\left(k_1,k_2;\frac{\alpha}b;z\right), \qquad z^{-k_1} F\left(k_1,k_1-\frac{\alpha}b+1; k_1-k_2+1,\frac1z\right),\] where $F(a,b;c;z)$ is the hypergeometric function. Since $k_1>\alpha/b>0$ and $\alpha/b>k_2$, one can show that both solutions are well-defined and linearly independent and thus form a fundamental system of solutions to the hypergeometric equation. We obtain the solution to  \eqref{newide}
\by \psi(x)=C_1 f_1(x)-C_2 f_2(x).\ey where 
\by
f_1(x)= F\left(k_1,k_2;\frac\alpha b; -\frac{c}{a} e^{bx}\right), \qquad f_2(x)=e^{-bk_1 x}F\left( k_1,k_1-\frac{\alpha}b+1; k_1-k_2+1;-\frac{a}{c}e^{-bx} \right),
\ey
subject to the initial condition that \be (a+c)\psi'_1(0)-(\lambda+\delta)\psi_1(0)=0.\label{bdrycond}\ee
It follows from \cite[p.387, (15.5.1)]{Olver:2010:NHMF}
\[\frac{\df}{\df z} F(a,b;c;z)=\frac{a b}{c} F(a+1,b+1;c+1;z)\]
that the derivative of $f_1$ is given by
\by f'_1(x)&=&\frac{\delta}{a} e^{bx} F\left(k_1+1,k_2+1;\frac \alpha b+1; -\frac c a e^{bx}\right);\\
f''_1(x)&=&\frac{b\delta}{a}e^{bx}F\left(k_1+1,k_2+1;\frac{\alpha}{b}+1;-\frac{c}{a}e^{bx}\right)\\
&&-\frac{\delta(\alpha bc+b^2c-\alpha \delta-b\lambda-b\delta)}{a^2(\alpha+b)}e^{2bx}F\left(k_1+2,k_2+2;\frac{\alpha}{b}+2;-\frac{c}{a}e^{bx}  \right).\ey
Using \cite[p.387, (15.5.3)]{Olver:2010:NHMF} with $n=1$,
\[\frac{\df}{\df z}  \big( (z^{a} F(a,b;c;z)\big)=a z^{a-1} F(a+1,b;c;z),\] we find the first two derivatives of $f_2$
\begin{displaymath}
\begin{aligned}
f'_2(x)=&\  -bk_1 e^{-bk_1 x} F\left(k_1+1,k_1-\frac{\alpha}{b}+1; k_1-k_2+1;-\frac a c e^{-bx} \right),\\
f''_2(x)=&\  (b k_1)^2 e^{-bk_1 x} F\left( k_1+1, k_1-\frac \alpha b+1; k_1-k_2+1; -\frac a c e^{-bx} \right)\\
&\ -\frac{ab^2 k_1(k_1+1)(k_1-\alpha/b+1)}{c(k_1-k_2+1)}e^{-b(k_1+1)x}F\left(k_1+2,k_1-\frac \alpha b+2;k_1-k_2+2; -\frac a c e^{-bx}  \right).
\end{aligned}
\end{displaymath}
Using \cite[p.388, (15.5.21)] {Olver:2010:NHMF}
\[c(1-z)\frac{\df }{\df z} F(a,b;c;z)=(c-a)(c-b)F(a,b;c+1;z)+c(a+b-c)F(a,b,c;z), \] we obtain
\be \label{bdrycomp1} (a+c)f'_1(0)-(\lambda+\delta)f_1(0)=-\lambda F\left(k_1,k_2; \frac \alpha b +1; - \frac c a \right).\ee
Using \cite[p.388, (15.5.13)]{Olver:2010:NHMF}
\by  (c-a-b)F(a,b;c;z)+a(1-z)F(a+1,b;c;z)-(c-b)F(a,b-1;c;z)=0,\ey we get
\be (a+c)f'_2(0)-(\lambda+\delta)f_2(0)=(bk_2-\alpha)c F\left( k_1,k_1-\frac \alpha b; k_1-k_2+1; -\frac a c\right). \label{bdrycomp2}\ee
Substituting \eqref{bdrycomp1}, \eqref{bdrycomp2} into \eqref{bdrycond} determines the unique solution to \eqref{newide} up to a multiplicative constant
\[\psi_1(x)=C_1 f_1(x)-C_2 f_2(x),\] where $C_1$ and $C_2$ are positive constants given by
\by
C_1 = \left(\alpha-bk_2\right)c F\left(k_1,k_1- \frac \alpha b; k_1-k_2+1;-\frac a c\right),\qquad C_2 = \lambda F\left( k_1, k_2; \frac \alpha b+1; -\frac c a\right).
\ey
We claim that a negative, increasing and concave solution to \eqref{newode} is given by
\[\psi_2(x)=- f^\ast_2(x):=- e^{-b k^\ast_1 x} F\left(k^\ast_1, k^\ast_1-\frac \alpha b +1; k^\ast_1-k^\ast_2+1; -\frac{a}{c-u_0} e^{-bx} \right),\] where
\by k^\ast_1=\frac{\alpha (c-u_0)-\lambda-\delta+\sqrt{(\alpha (c-u_0)-\lambda-\delta)^2+4\alpha\delta (c-u_0)}}{2b(c-u_0)}, \\
 k^\ast_2=\frac{\alpha (c-u_0)-\lambda-\delta-\sqrt{(\alpha (c-u_0)-\lambda-\delta)^2+4\alpha\delta (c-u_0)}}{2b(c-u_0)}.\ey
It follows from the definition of the hypergeometric function that $f^\ast_2(x)>0$ and $(f^\ast_2)'(x)<0$ for large $x$. Suppose there exists $x_0:=\inf\{x\ge 0: (f^\ast_2)'(x)=0, f^\ast_2(x)\ge 0\}.$ Recall that $f^\ast_2$ is a solution to the ODE \eqref{newode}. Since $a e^{-bx_0}+c-u_0>0$ and $-\alpha \delta<0,$ we must have $(f^\ast_2)''(x_0)\ge 0$ which leads to a contradiction. Thus $(f^\ast_2)'(x)<0, f^\ast_2(x)>0$ for all $x\ge 0$. To show the sign of $(f^\ast_2)''$, we differentiate \eqref{newode} with respect to $x$ and observe that $f^\ast_2$ satisfies
\[A(x) f'''(x)+B(x)f''(x)+C(x)f'(x)=0,\]
where $A(x)=a e^{-bx}+c-u_0>0, C(x)=ab(b-\alpha)e^{-bx}-\alpha \delta<0$ for $x\ge 0$ due to \eqref{newcond2}. We can conclude in the same way as before that $(f^\ast_2)''(x)>0$ for all $x\ge 0$. Therefore, $\psi_2$ is indeed negative, increasing and concave.

\subsubsection{Restricted Payment Scheme}

Consider the restricted payment scheme under which the dividend rate is capped at $0<u_0<c$. Then we can determine the optimal dividend strategy according to Theorem \ref{finalthr}. The condition \eqref{unique} reduces to
\be \frac{F(k^\ast_1, k^\ast_1-\alpha/b+1;k^\ast_1-k^\ast_2+1;-a/(c-u_0))}{bk^\ast_1 F(k^\ast_1+1,k^\ast_1-\alpha/b+1;k^\ast_1-k^\ast_2+1;-a/(c-u_0))}
< \frac{u_0(\lambda+\delta)-(a+c)\delta}{(\lambda+\delta)\delta}.\label{dexist}\ee
Therefore, under conditions \eqref{newcond1}, \eqref{newcond2}, \eqref{dexist}, the optimal value function is given by
\by 
V_R(x)=\left\{ \begin{array}{ll} \displaystyle \frac{C_1 f_1(x)-C_2 f_2(x)}{C_1 e_1(d)+C_2 e_2(d)}, & \mbox{ if } 0 \le x <d;\\
\displaystyle \frac{u_0}{\delta}-\frac{f^\ast_2} {e_3(d)}, & \mbox{ if } x \ge d,\end{array} \right.
\ey where $e_1, e_2, e_3$ are positive functions given by
\by
e_1(d)&=& \frac \delta a e^{bd} F\left(k_1+1, k_2+1; \frac{\alpha}{b}+1; -\frac{c}{a}e^{bd}   \right);\\
e_2(d)&=& bk_1 e^{-bk_1 d} F\left(k_1+1,k_1-\frac{\alpha}{b}+1; k_1-k_2+1;-\frac a c e^{-bd} \right);\\
e_3(d)&=& bk^\ast_1 e^{-bk_1 d} F\left(k^\ast_1+1,k^\ast_1-\frac{\alpha}{b}+1; k^\ast_1-k^\ast_2+1;-\frac{a}{c-u_0} e^{-bd} \right).
\ey
If \eqref{dexist} is not true, then the value function is given by
$$V_R(x)=\frac{u_0}{d}-\frac{f^\ast_2(x)}{K}, \qquad x\ge 0, $$
where 
\by
K=\frac{\delta  (b k^\ast_2-\alpha)(c-u_0)}{\lambda u_0}F\left(k^\ast_1, k^\ast_1-\frac{\alpha}b; k^\ast_1 -k^\ast_2 +1;-\frac{a}{c-u_0}\right).
\ey

\subsubsection{Unrestricted Payment Scheme}
In this case there is no restriction on the dividend payment rate. According to Theorem \ref{thm-soln-qvi}, we first determine the minimum of $\psi'_1$ at $x^\ast$ by
 \begin{displaymath}
 \begin{aligned}
\lefteqn{C_1\frac{b\delta}{a}e^{bx^\ast}F\left(k_1+1,k_2+1;\frac{\alpha}{b}+1;-\frac{c}{a}e^{bx}\right)}\\
 &-C_1\frac{\delta(\alpha bc+b^2c-\alpha \delta-b\lambda-b\delta)}{a^2(\alpha+b)}F\left(k_1+2,k_2+2;\frac{\alpha}{b}+2;-\frac{c}{a}e^{bx^\ast}  \right)\\
&-C_2(b k_1)^2 e^{-bk_1 x^\ast} F\left( k_1+1, k_1-\frac \alpha b+1; k_1-k_2+1; -\frac a c e^{-bx^\ast} \right)\\
&+C_2\frac{ab^2 k_1(k_1+1)(k_1-\alpha/b+1)}{c(k_1-k_2+1)}e^{-b(k_1+1)x^\ast}F\left(k_1+2,k_1-\frac \alpha b+2;k_1-k_2+2; -\frac a c e^{-bx^\ast}  \right)=0.
\end{aligned}\end{displaymath}
The existence and uniqueness of such a value $x^\ast$ is guaranteed by conditions \eqref{newcond1} and \eqref{newcond2} according to Theorem \ref{finalbar}(i). Thus the optimal value function is given by
\by 
V(x)=\left\{ \begin{array}{ll} \displaystyle \frac{C_1 f_1(x)-C_2 f_2(x)}{C_1 e_1(x^\ast)+C_2 e_2(x^\ast)}, & \mbox{ if } 0 \le x <x^\ast;\\
\displaystyle x-x^\ast-\frac{C_1 f_1(x^\ast)-C_2 f_2(x^\ast)}{C_1 e_1(x^\ast)+C_2 e_2(x^\ast)}, & \mbox{ if } x \ge x^\ast,\end{array} \right.
\ey
In the case where $\alpha \lambda(a+c)-(ab+\lambda+\delta)(\lambda+\delta)\ge 0$, Theorem \ref{finalbar}(ii) tells us that the barrier level is $0$ and the optimal value function is given by
\[V(x)=x+\frac{a+c}{\lambda+\delta},\qquad x \ge 0.\]

\begin{figure}[htbp]
\begin{center}
 \mbox{\subfigure[Value functions $V(x)$ and $V_R(x)$]{\epsfig{figure=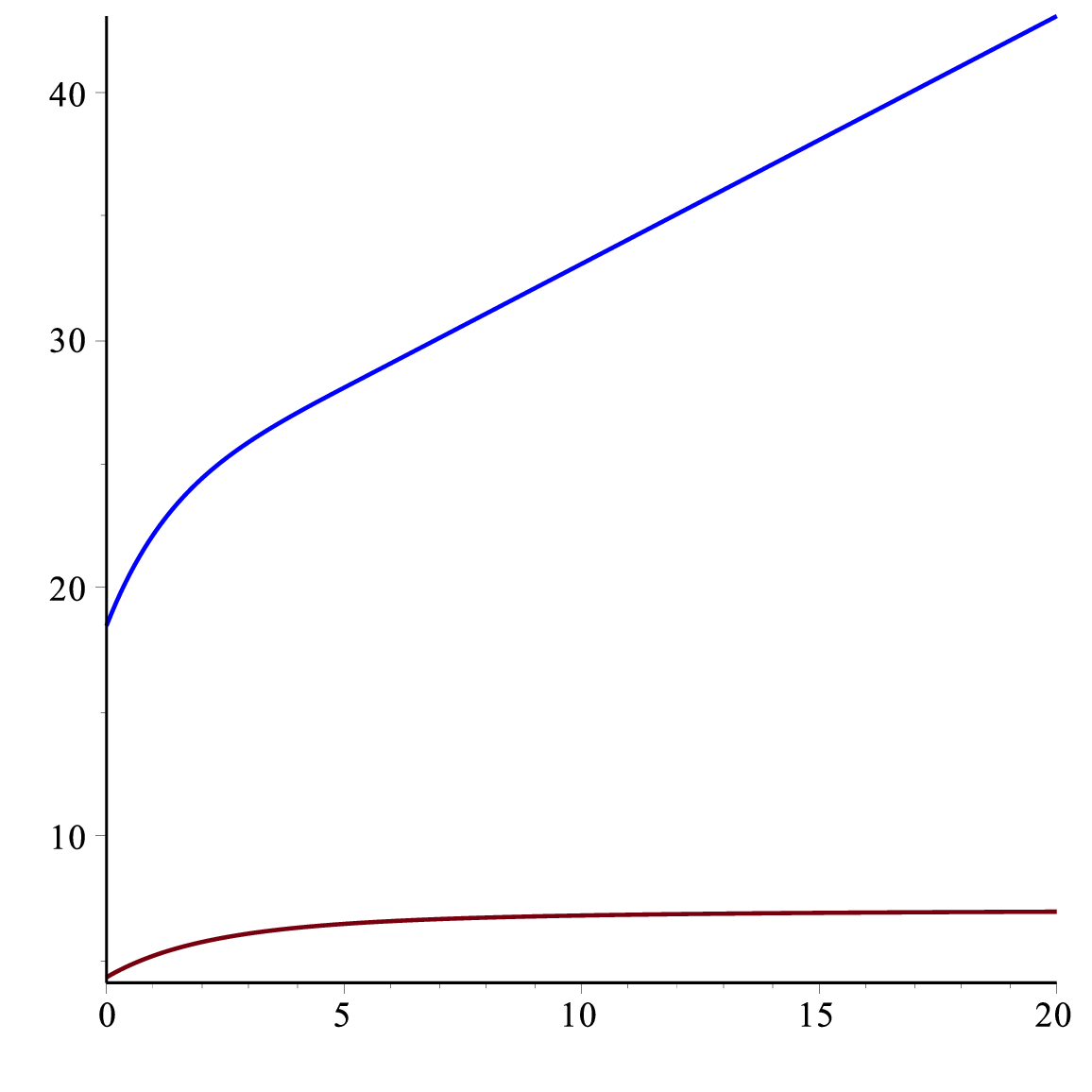,width=0.43\linewidth}}}\ \
 \mbox{\subfigure[The difference $V(x)-V_R(x)$]{\epsfig{figure=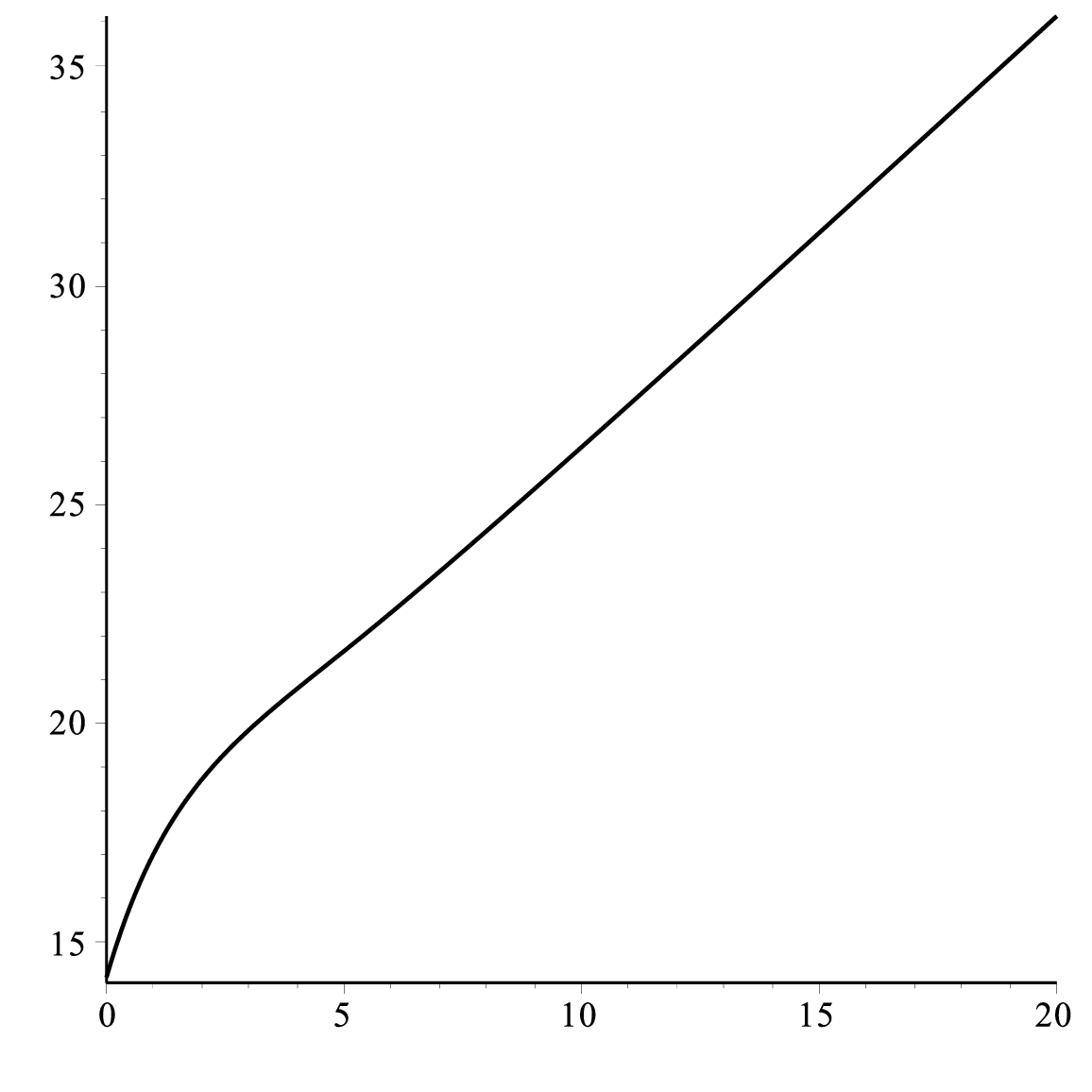,width=0.43\linewidth}}}
\caption{Comparison %of optimal policies
 in the regressive growth model} \label{fig3}
\end{center}
\end{figure}

Here we give a numerical example of the optimal dividend payment schemes with regressive premium rates. The parameters are chosen as follows: $a= 2, b=0.15, c= 2, \lambda= 1, \alpha= 1, $ and $\delta=0.06$. All conditions \eqref{newcond1}, \eqref{newcond2}, \eqref{dexist} are satisfied. The values of $d=0.2143382330$ and $x^\ast=5.324276475$ are determined numerically by a root search method. The value functions $V(x),  V_R(x)$, and their difference are plotted in Figure \ref{fig3}.

\section{Conclusions and Remarks}\label{sect-conclusion}

This work is devoted to the optimal dividend payment problems
for the %piecewise-deterministic compound Poisson 
PDCP risk model. 
Both restricted and unrestricted dividend schemes are investigated and compared. 
We provide easily verifiable conditions under which the threshold and barrier strategies are optimal 
restricted and unrestricted dividend payment policies, respectively. 
 % As with many well-known results in other diffusion risk models,
% We prove a common optimality which has been observed previously in many special cases
 % in the setting of PDCP risk models
% that the optimal restricted dividend payment scheme is the threshold strategy,
% in which dividends are paid only at the maximal rate as long as the surplus remains above a certain level; and the optimal unrestricted dividend payment scheme is the barrier strategy, in which dividends are only paid at times when the surplus reaches a (possibly different)  level and at such a rate that the surplus stays at the same level until the next claim arrival or ruin.
%Rather than treating different special cases separately,
 Our analysis is primarily based on the qualitative properties of solutions to certain IDE and ODE associated with a general PDCP risk model. Three examples are demonstrated to illustrate the main results. 

A number of   questions deserve further investigations. % in future work.
One can consider more general models in which the parameters and hence the dynamics of the
surplus level depend on a stochastic process such as a continuous-time Markov chain. In such a case, we need to deal with regime-switching
jump diffusions (see \cite{YZ-10}) and the resulting HJB equation will be a coupled system of nonlinear integro-differential  equations. It is conceivable that
it will be more challenging to
obtain
the corresponding value function and optimal dividend policy in closed forms.
Some initial work in this vein can be found in \cite{Zhu-10} using the viscosity solution approach.
% It will also be interesting to approach the optimal dividend payment problem in the setting of piecewise deterministic compound Poisson risk model
 % using the powerful viscosity solution framework (Crandel et al. \cite{CIL92}).
Another problem of interest is to consider transaction costs, reinsurance, and/or  investments.
Similar work in the setting of controlled diffusions can be found in \cite{Bai-Paulsen-10,Choulli}, etc.

%\section*{Acknowledgement}
% Chao Zhu is grateful for financial support from the National Science Foundation
% under grant DMS-1108782 and from the UWM Research Growth Initiative, and City University
% of Hong Kong (SRG) 7002677.

\appendix
\section{Appendix}%{Qualitative Analysis}

We prove two technical lemmas on the properties of the solution to the ODE \eqref{sode}, which will enable us to identify optimal dividend strategies in a general PDCP model.
Consider
\begin{equation}\label{sode} f(x) u''(x)+h(x) u'(x)-k u(x)=0,\qquad x \ge 0\end{equation} where $f,h:[0,\infty)\rightarrow \rr$ are continuously differentiable and $f(x) >0$ for all $x \ge 0$ and $k>0$.
 % We identify two sets of sufficient conditions under which solutions can be obtained for optimal dividend policies.

\begin{lem} \label{idecond}
Let $u$ be a solution of \eqref{sode} such that
\begin{itemize}
    \item[(i)] $u(0)>0, u'(0)>0, u''(0)<0;$
%    \item There are constants $0<A<1$ and $B>0$ such that $h(x)/k \le Ax+B$ for all $x \ge 0$;
    \item[(ii)] There exists $m<k$ such that $h'(x) \le m $ for all $x>0$.
\end{itemize}
Then there is a unique $b>0$ such that $u(x)$ is concave for $0 \le x<b$ and convex for $x>b$.
\end{lem}
\begin{proof} 
%Note that a similar result was proved in \cite[Lemmas 4.1 and 4.2]{Shreve-84}. But their result can not be directly applied in our setting. Therefore we include a brief proof here, which is different from their proof and will facilitate our later presentation. 

 First, we show that $u(x)>0$ and $u'(x)>0$ for all $x\ge 0$. Since $u(0)>0$ and $u'(0)>0,$  $u(x)>0$ and $u'(x)>0$ for small $x>0.$ Suppose there exists $x_0>0$ such that $u'(x_0)=0$ and $u'(x)>0$ for all $x\in (0,x_0)$, then $u''(x_0)\le 0$. It follows from \eqref{sode} that $u''(x_0)=ku(x_0)/f(x_0)>0$ which leads to a contradiction. Hence the claim is proved.

Then, we prove by contradiction that there exists $b>0$ such that $u''(b)=0$. Suppose that $u''(x)$ has no zero. Since $u''(0)<0$,  $u''(x)<0$ for all $x \ge 0$. Let $A=\max\{m/k, 1/2\}<1$. Since $h'(x) \le m$ for all $x>0$,   there must exist $B>0$ such that $h(x)/k \le (m/k)x+h(0)/k \le Ax+B.$  It follows from \eqref{sode} that
\begin{equation}(Ax+B)u'(x) \ge \frac{h(x)}{k} u'(x) > u(x) \qquad \mbox{ for all } x \ge 0.\label{ine}\end{equation}
Dividing \eqref{ine} by $(Ax+B)u(x)>0$ and integrating from $0$ to $x$, we obtain
$ u(x) \ge u(0) \left( \frac{Ax+B}{B} \right)^{1/A},$ which implies
\[u'(x) \ge (Ax+B)^{-1}u(x) \ge u(0) B^{-1/A}(Ax+B)^{1/A-1} \rightarrow \infty, \qquad x \rightarrow \infty.\] This is impossible because $u'(x)$ is a decreasing function. The contradiction shows that $u''(x)$ has a positive zero.

Last, we show that $u''(x)$ has a unique zero. Suppose that $u''(x)$ does have a second zero $x_1>b$ such that $u''(x)>0$ for all $x\in (b,x_1)$ and $u''(x_1)=0$. Hence it must be true that $u'''(x_1)\le 0$. Differentiating \eqref{sode}, we obtain
\begin{equation} f(x)u'''(x)+(f'(x)+h(x))u''(x)+(h'(x)-k)u'(x)=0.\label{tode}\end{equation} Therefore, we must also have $u'''(x_1)=(k-h'(x_1))u'(x_1)/f(x_1)>0$ since $u'(x)>0$ for all $x \ge 0$. This contradiction implies the uniqueness of the zero of $u''(x)$.
\end{proof}

\begin{lem} \label{odecond}
Assume that there exists a constant $m<k$ such that $h'(x) \le m$ for all $x \ge 0$. Then equation \eqref{sode} admits a solution $u$ such that
\[u(x)<0, \, u'(x)>0, \, u''(x)< 0 \,\mbox{ for all } x \ge 0.\]
\end{lem}
\begin{proof} Let $v(x)=u'(x)$ in equation \eqref{tode}. We obtain
\beq f(x)v''(x)+(f'(x)+h(x))v'(x)+(h'(x)-k)v(x)=0.\label{vode}\eeq
Since $h'(x)\le m<k$ for all $x\ge 0$, it follows from Corollary 1.2 in Chapter XIV of \cite{Hartman} that \eqref{vode} has at least one solution $v(x)$ satisfying $v(x)>0, v'(x) \le 0$ for all $x \ge 0$. If there exists $x_2 \ge 0 $ such that $v'(x_2)=0$, then we must have $v''(x_2)\le 0$. Note that it follows from \eqref{vode} that $v''(x_2)=(k-h'(x_2))v(x_2)/f(x_2)>0,$ which causes a contradiction. Therefore, we must have $v(x)>0, v'(x)<0$ for all $x \ge 0$. There is a function $u(x)$ such that $u'(x)=v(x)$ and $u(x)$ satisfies \eqref{sode}. Thus, for such a solution $u(x)$, we must have $u'(x)>0, u''(x)<0$ for all $x \ge 0$.

Next, we show that $u(x)<0$ for all $x \ge 0$. Suppose there exists $x_3 \ge 0$ such that $u(x_3) \ge 0$. Then we must have $u(x_4)>0, u'(x_4)>0, u''(x_4)<0$ for some $x_4>x_3$. Let $w(x)=u(x-x_4), f_w(x)=f(x-x_4), h_w(x)=h(x-x_4).$ Then $w(x)$ must satisfy the ODE
 \[f_w(x)w''(x)+h_w(x) w'(x)-kw(x)=0, \qquad x \ge 0,\] and all conditions in Lemma \ref{idecond} are satisfied. Hence there must exist $x_5>0$ such that $u''(x_4+x_5)=w''(x_5)=0$ which contradicts the fact that $u''(x)<0$ for all $x \ge 0$. Therefore, it must be true that $u(x)<0$ for all $x \ge 0$.
\end{proof}

\begin{proof}[{Proof of Lemma \ref{lem-bounds}}]
 That the function $V_R(x)$  is bounded by $u_{0} / \delta $ is obvious. Also, monotonicity can be established by considering two surplus processes with different initial surplus levels and the same dividend payment schemes. 
The rest of the proof is divided into two steps.

% {\em Step 1.} (Monotonicity) Let $x_2 > x_1 \ge 0$. 
% Denote by $X_1(t)$ the surplus process
%  under the dividend payment scheme
 %  $D_1=\{\int^t_0 u_1(s)\df s, t\geq 0\}$ and initial surplus level  $x_1$.
 % Put $ \tau_1:= \inf\set{t\ge 0: X_1(t) < 0}. $
 % Also, let $D_2=\{\int^t_0 u_2(s)\df s, t\geq 0\}$ be such that $u_2(t) =u_1(t)$ for all $t\le \tau_1$.
 % Let $X_2(t)$ be the surplus process
 % under the dividend payment scheme $u_2\cd$ and initial surplus level  $x_2$  and $ \tau_2:= \inf\set{t\ge 0: X_2(t) < 0}. $
 % Then apparently we have $\tau_2 \ge \tau_1$.
 % Consequently it follows that
 % $ V_R(x_2) \ge J(x_2,D_2) =\ex \int_0^{\tau_2} e^{-\dl t} u_2(t) \df t \ge \ex\int_0^{\tau_1} e^{-\dl t} u_1(t) \df t = J(x_1,D_1).$
 % Finally, taking superemum over $ D_1 \in \Pi_R$, it follows that $V_R(x_2) \ge V_R(x_1)$, as desired.

{\em Step 1.} (Lipschitz Continuity) Let $h$ be a small positive value and
  % $\td{\pi}:=\{\int^t_0\td{u}(s)\df s, t\geq 0\}$
$\tilde D \in \Pi_R$ be an arbitrary strategy. Define another strategy
  % $\pi:=\{\int^t_0 u(s)\df s, t\geq 0\}$
$D\in \Pi_R $ so that $ D(t) = \td D(t-h)I_{\set{T_1 \wedge t > h}}.$
 % \[D(t)=\left\{ \begin{array}{ll} \td{D}(t-h), &\qquad \mbox{if } 
 %  T_1 \wedge t > h ;\\ 0, & \qquad \mbox{otherwise.}\end{array} \right.\]
Denote  by $X(t) $ the surplus process
with initial surplus $x>0$ under the dividend payment strategy $D$.
 It is clear that
if $T_1 >h$, then the surplus at time $h$ is $X(h) =\phi_x(h)$.
Therefore it follows that
$V_R(x) \ge J(x,D) \geq e^{-(\lambda+\delta)h} J(\phi_x(h), \td D),$ which implies by taking the supremum over all possible strategies $\td D\in \Pi_R$ that
\begin{eqnarray} V_R(x) \geq e^{-(\lambda+\delta)h} V_R(\phi_x(h)) \geq e^{-(\lambda+\delta)h} V_R(x),\label{inter}\end{eqnarray} with the last inequality from the fact that $V_R(x)$ is an increasing function and $g(x)\geq 0.$ Thus $V_R(x)$ is right continuous by the continuity of $\phi_x(h)$ and
the squeeze theorem. Replacing $x$ with $\phi_x(-h)$ in \eqref{inter}, we  obtain
\begin{eqnarray*} V_R(\phi_x(-h)) \geq e^{-(\lambda+\delta)h} V_R(x) \geq e^{-(\lambda+\delta)h} V_R(\phi_x(-h)).\end{eqnarray*} Hence left continuity follows.
 Now it follows from  \eqref{inter} that
 \[0\le V_R(\phi_x(h))-V_R(x)\leq (1-e^{-(\lambda+\delta)h}) V_R(\phi_x(h))
 \leq  (1-e^{-(\lambda+\delta)h}) u_0/\delta.\] Therefore, $V_R(x)$ is indeed Lipschitz continuous.

{\em Step 2.} (Limit at $\infty$) Let $D(t):= u_0 t$ for all $t\ge 0$. Denote the
surplus process by $X(t)$ under the strategy $D$ and initial surplus $x>0$ and by $\tau$ the corresponding time of ruin. Then as $x\to \infty$, $\tau$   converges to infinity a.s.
Therefore
\bed V_R(x) \ge J(x,D) = \ex\int_0^\tau e^{-\dl t} u_0 \df t = \frac{u_0}{\dl} (1-\ex[e^{-\dl \tau}]) \to \frac{u_0}{\dl} ,\eed this, together with the boundedness of $V_R(x)$, leads to the desired conclusion.
\end{proof}

\begin{proof}[Proof of Theorem \ref{thm-HJB}] The proof is  motivated by \cite[Theorem 2.32]{schmidli2008}.
  It utilizes Lemma \ref{lem-bounds}  % some renewal type argument, 
   and the dynamic programming principle: 
\begin{eqnarray}\label{DPP}
 V_R (x)=\sup_{u\cd\in \Pi_R}
 \mathrm{E}_{x}\left[\int_0^{ \tau \wedge \theta } e^{-\delta s}
  u(s) \df s +  e^{-\delta (  \theta \wedge  \tau)}V_R(X (\theta \wedge \tau  ) )   \right],\ \ x\geq 0,
\end{eqnarray}  where $\theta$ is an $ \mathcal {F}_{t}$-stopping time.

{\em Step 1.} Let $h>0$ and   $u\in [0,u_0]$.
Let $\{ D(t):=ut, t\geq 0\} \in \Pi_R.$
Denote $$\phi(t,x)= x + \int_0^t (g(\phi(s,x))-u)\df s, \  t\ge 0.$$
 % {\red $\tilde{\pi}:= \{\int^t_0 \tilde{u} (X(s)) \df s, t\geq 0\}$
 % be an arbitrary strategy}, and
 % $ {\red \tilde{\phi}} (h,x) =  x+\int^{h}_{0}(g(X_{t})-  \tilde{u}
 % (X_{t}))d t\geq
 % 0.      % $
 % \begin{eqnarray}\nonumber
 % \phi (h,x) =  x+\int^{h}_{0}(g(X_{t})-  \tilde{u} (X_{t}))d t\geq
 %  0.
 % \end{eqnarray}
 % Choose $\varepsilon>0$
 %  We consider the following strategy:
 % \begin{eqnarray}\nonumber
 % r(X_{t})=\left\{\begin{array}{ll}
 %  \tilde{u}(X_{t}), ~~~~~~~~~ 0\leq t\leq\tau_{1}\wedge h,\\
 % r^{\varepsilon}(X_{t-\tau_{1}\wedge h}), ~~ t>\tau_{1}\wedge h.
 % \end{array}
 % \right.
 %\end{eqnarray}
 % $\{r^{\varepsilon}(X_{t})\}$ is a strategy for initial capital
 % $x_{k}$, $x_{k}\leq X_{\tau_{1}\wedge h}<x_{k+1}$, such that
 % $V_R^{\varepsilon}(x_{k})>V_R(x_{k})-\varepsilon$. Where
 % $x_{k}=\frac{k}{n}\bar{\phi}(h,x)$ for $0\leq k\leq n$.
 %  Thus,
 % $\{r(X_{t})\}$ is measurable.
 % By the Lipschitz continuity of $V_R(x)$, we can choose $n$ large
 % enough such that
 % $V_R^{\varepsilon}(\tilde{x})>V_R(\tilde{x})-\varepsilon$ for all
 % $\tilde{x}\in [0,\bar{\phi}(h,x)]$.
 The arrival time $T_1$ of the first claim 
  % has density $\lambda e^{-\lambda t}$ and $
 % T_{1}$ is larger than $h$ with probability $ e^{-\lambda h}$. 
is exponentially distributed with parameter $\lambda$.  Hence we can use
the law of total probability and take $\theta= T_{1}\wedge h$ in
\eqref{DPP} to obtain
\beq\label{DPP-1}\begin{aligned}
V_R(x)
  \geq &\  e^{-\lambda h}\bigg[\int^{h}_{0}e^{-\delta
t}u \df t+e^{-\delta h}V_R ( \phi
(h,x))\bigg]  \\ &  +\int^{h}_{0}\lambda e^{-\lambda
t} \bigg[\int^{t}_{0}e^{-\delta
s}u \df s +e^{-\delta t}\int^{
\phi (t,x)}_{0}V_R (  \phi (t,x)-y)\df Q(y)\bigg] \df
 t.
\end{aligned}\eeq
Rearranging the terms and dividing by $h$ yields
\begin{equation}\label{DPP-S}
\begin{aligned}
  & \frac{V_R( \phi
(h,x))-V_R(x)}{h}-\frac{1-e^{-(\lambda+\delta)h}}{h}V_R(  \phi
(h,x))+\frac{e^{-\lambda
h}}{h}\int^{h}_{0}e^{-\delta t}u\df  t\\
& \  +\frac{1}{h}\int^{h}_{0}\lambda e^{-\lambda
t}\bigg[\int^{t}_{0}e^{-\delta
s}u\df  s+e^{-\delta t}\int^{ { \phi}
(t,x)}_{0}V_R( \phi (t,x)-y)\df Q(y)\bigg]\df  t \leq
0. \end{aligned}
\end{equation}

Let
$$ \displaystyle
D^{+}{V_R}(x)= \limsup_{ \Delta \rightarrow
0+}\frac{V_R(x+ \Delta)-V_R(x)}{ \Delta },
\,\,\, D^{-} V_R (x)=\liminf_{ \Delta \rightarrow 0+}\frac{V_R(x+
\Delta  )-V_R(x)}{ \Delta }.
$$
Note that $  D^{+}{V_R}(x) $ and   $  D^{-}{V_R}(x)  $  are finite by the Lipschitz
continuity of $V_R$. Recall the facts that $0\le u\le u_0 < g(x)$,  $g$ is continuous,
$\phi(\cdot, x)$ is strictly increasing and $\phi(t,x)\to x$ as $t\to 0$. Hence we can  choose a sequence $\set{h_n, n\ge 1}$ satisfying $h_n \to 0 $ as $n \to \infty$ and
\bed \lim_{n\to \infty} \frac{V_R(\phi(h_n,x))-V_R(x)}{\phi(h_n,x)-x} = D^+ V_R(x). \eed
By the definition of $\phi(h_n,x)$, we have $\phi(h_n,x) \to x$ as $n\to \infty$. Also, we have from the continuity of $g$ that as $n\to \infty$
\bed \frac{\phi(h_n,x) - x}{h_n} =\frac{1}{h_n}\int_0^{h_n} \left[g(X(t))-u\right] \df t \to g(x)-u. \eed
Hence  it follows that
\beq\label{D+V(x)} \begin{aligned}\lim_{n\to \infty} \frac{V_R(\phi(h_n,x))-V_R(x)}{h_n} & =
 \lim_{n\to \infty} \frac{V_R(\phi(h_n,x))-V_R(x)}{\phi(h_n,x) - x}\cdot\frac{\phi(h_n,x) - x}{h_n}= D^+V_R(x)(g(x)-u).
\end{aligned} \eeq
 Now taking $h=h_n$ in \eqref{DPP-S} and letting $n\to \infty$, in view of \eqref{D+V(x)},
 detailed calculations reveal  that
 \beq \label{eq-HJB-1}
 [g(x)-u]D^+ V_R(x) - (\lambda+\delta) V_R(x) + \lambda \int_0^x V_R(x-y)\df Q(y) + u \le 0, \ \ \forall u \in [0,u_0].
 \eeq

{\em Step 2.} On the other hand, by the definition of
$V_R$ in \eqref{restricted-V-defn}, there exists a strategy $\bar{D}:=\{\int^t_0 \bar{u}(s)\df s, t\geq 0 \}\in \Pi_R$ such that
$ J(x,\bar D)\geq V_R(x)-h^{2}$.
 Denote
 $$
 \bar{\phi } (t,x)= x+\int^{t}_{0}(g(X(s))-\bar{u} (s))\df  s. $$
 Then as argued before, we must have
 \bed \begin{aligned}
 V_R(x) & \le J(x,\bar D)+ h^2 \\
  & \le h^2 + e^{-\la h} \left[\int_0^h e^{-\delta s} \bar u(s)\df s + e^{-\delta h} V_R(\bar \phi(h,x))\right] \\
  & \ \ +   \int^{h}_{0}\lambda e^{-\lambda
t} \bigg[\int^{t}_{0}e^{-\delta
s}\bar u(s) \df s +e^{-\delta t}\int^{\bar \phi (t,x)}_{0}V_R ( \bar\phi (t,x)-y)\df Q(y)\bigg] \df t.
 \end{aligned}\eed
 We find by rearranging the terms and dividing by $h$ that
\begin{equation} \label{DPP-L}
\barray
 \ad h+ \frac{V_R( \bar{\phi}
 (h,x))-V_R(x)}{h}-\frac{1-e^{-(\lambda+\delta)h}}{h}V_R( \bar{ \phi  }  (h,x)) +\frac{e^{-\lambda
h}}{h}\int^{h}_{0}e^{-\delta t} \bar{u}  (t)\df  t \\ \aad
 +\frac{1}{h}\int^{h}_{0}\lambda
e^{-\lambda t}\bigg[\int^{t}_{0}e^{-\delta s}\bar{u}
(s)\df s+e^{-\delta t}\int^{
     \bar{\phi  }     (t,x)}_{0}V_R(\bar{\phi}    (t,x)-y)\df Q(y)\bigg]\df  t\geq 0. \earray
\end{equation}
Denote
$\bar u  := \liminf_{s\to 0+} \bar u(s) \in [0,u_0].$
Then it follows %from Fatou's Lemma
that
\bed \limsup_{h\to 0+} \frac{\bar \phi(h,x)-x}{h} = \limsup_{h\to 0+} \frac{1}{h } \int_0^h [g(X(s))-\bar u(s)]\df s  \le g(x)-\bar u.\eed
As in Step 1, we can choose a sequence $h_{m} \rightarrow 0+$ such that
\begin{eqnarray}\nonumber
\lim_{m\rightarrow\infty}\frac{V_R( \bar{\phi}( h _{m},x))-V_R(x)}{
\bar{\phi}( h_{m},x)-x      }= D^{-}V_R(x ).
\end{eqnarray}
Note that $D^- V_R(x) \ge 0$ by the monotonicity of $V_R$. Then it follows that
\bed \barray \disp\limsup_{m\to \infty}
  \frac{V_R(\bar{\phi}(h_{m},x))-V_R(x)}{h_m} \ad =\limsup_{m\to \infty}
  \frac{V_R(\bar{\phi}(h_{m},x))-V_R(x)}{\bar\phi(h_{m},x)-x}\cdot \frac{\bar\phi(h_{m},x)-x}{h_m} 
  \\ \ad \le D^-V_R(x) [g(x)-\bar u].\earray
\eed
Now  by taking $h=h_m$ and letting $m\to \infty$ in \eqref{DPP-L}, we obtain
\beq\label{eq-HJB-2}
[g(x)-\bar u] D^-V_R(x) - (\la +\dl) V_R(x) + \la \int_0^x V_R(x-y)\df Q(y) + \bar u \ge 0.
\eeq

{\em Step 3.}
Note that $g(x)- \bar u \ge g(x) -u_0 > 0$.
 % $$ D^-V_R(x) \ge \frac{(\la +\dl) V_R(x) - \la \int_0^x V_R(x-y)\df Q(y)
 % - \bar u}{g(x)-\bar u}.$$
 Hence,   by taking $u=\bar u $ in
\eqref{eq-HJB-1} and comparing the resulting equation
with \eqref{eq-HJB-2}, we have
$D^+ V_R(x) \le D^- V_R(x).$
But
$D^+V_R(x) \ge D^-V_R(x)$ by definition. Thus it follows that
$D^+V_R(x) = D^-V_R(x)$ or $V_R(x)$ is differentiable from the right. Moreover, a combination of \eqref{eq-HJB-1} and \eqref{eq-HJB-2}
yields that $V_R'(x+)$, the right derivative of $V_R$, satisfies the HJB equation
\begin{equation}\label{HJB-right}
 \sup_{0\leq u  \leq u_{0} }\bigg\{[g(x)
-u ]V_R'(x+)-(\lambda+\delta)V_R(x) +\lambda\int^{x}_{0}V_R(x-y)\df
Q(y)+u \bigg\}=0.
\end{equation}

Similarly, we obtain that the left derivative $V_R'(x-)$ exists and fulfills the HJB
equation
\begin{equation}\label{HJB-left}
 \sup_{0\leq u  \leq u_{0}}\set{[g(x)
-u]V_R'(x-)-(\lambda+\delta)V_R(x)+\lambda\int^{x}_{0}V_R(x-y)\df
Q(y)+u}=0.
\end{equation}

{\em Step 4.}
With \eqref{HJB-right} in hands, we claim that
\begin{equation}\label{claim-1}  V_R'(x+)\gtreqqless 1
\Leftrightarrow
 %\text{ if and only if }
(\lambda+\delta)V_R(x)-\lambda\int^{x}_{0}V_R(x-y)\df
Q(y)\gtreqqless g(x).
\end{equation}
In fact, if $V_R'(x+) > 1$, then
\bea 0\ad = \sup_{0\le u \le u_0} \set{[g(x)
-u ]V_R'(x+)-(\lambda+\delta)V_R(x) +\lambda\int^{x}_{0}V_R(x-y)\df
Q(y)+u  }\\
\ad = g(x) V_R'(x+) -  (\lambda+\delta)V_R(x) +\lambda\int^{x}_{0}V_R(x-y)\df
Q(y).
\eea
Hence we have
$ (\lambda+\delta)V_R(x) - \lambda\int^{x}_{0}V_R(x-y)\df Q(y) = g(x) V_R'(x+) > g(x).$
Conversely, if $ (\lambda+\delta)V_R(x) - \lambda\int^{x}_{0}V_R(x-y)\df Q(y) > g(x)$, then we have
\bea 0\ad = \sup_{0\le u \le u_0} \set{[g(x)
-u ]V_R'(x+)-(\lambda+\delta)V_R(x) +\lambda\int^{x}_{0}V_R(x-y)\df
Q(y)+u  } \\
\ad < \sup_{0\le u \le u_0} \set{(g(x)-u)(V_R'(x+)-1)}.\eea
But $g(x)-u \ge g(x) -u_0 >0$. Thus we must have $V_R'(x+) >1$.
Hence the first case in \eqref{claim-1} follows.
Similar arguments establish the other two cases in \eqref{claim-1}.

Similarly, \eqref{HJB-left} leads to
\begin{equation}\label{claim-2}
V_R'(x-)\gtreqqless 1 \Leftrightarrow
(\lambda+\delta)V_R(x)-\lambda\int^{x}_{0}V_R(x-y)\df
Q(y)\gtreqqless g(x) .
\end{equation}

Hence, \eqref{claim-1} and \eqref{claim-2} imply that $V_R'(x+)$ and $V_R'(x-)$ are both less than $1$, both greater
than $1$, or both equal to $1$.
This, together with the HJB equations
(\ref{HJB-right}) and (\ref{HJB-left}),
implies  that $V_R'(x+)=V_R'(x-)$ and so $V_R'(x)$ exists. Moreover, the continuities of $V_R$ and $g$ implies that $V_R'(x)$ is continuous. That is, $V_R(x)$ is continuously
differentiable and satisfies the HJB equation (\ref{HJB}).

{\em Step 5.} Finally, the optimality of the strategy \eqref{eq:u_R} is obvious since \eqref{HJB} is linear in $u$. 
\end{proof}

%\bibliography{/Users/zhu/papers/refs}\end{document}

\begin{thebibliography}{}
\parskip=-3pt

\bibitem[Albrecher and Hartinger, 2007]{Alb-Har-2007}
Albrecher, H. and Hartinger, J. (2007).
\newblock A risk model with multilayer dividend strategy.
\newblock {\em N. Am. Actuar. J.}, 11(2):43--64.

\bibitem[Albrecher and Thonhauser, 2008]{Albrecher-Thon-08}
Albrecher, H. and Thonhauser, S. (2008).
\newblock Optimal dividend strategies for a risk process under force of
  interest.
\newblock {\em Insurance Math. Econom.}, 43(1):134--149.

\bibitem[Asmussen et~al., 2000]{Asmussen-H-T-00}
Asmussen, S., H{\o}jgaard, B., and Taksar, M. (2000).
\newblock Optimal risk control and dividend distribution policies. {E}xample of
  excess-of loss reinsurance for an insurance corporation.
\newblock {\em Finance Stoch.}, 4(3):299--324.

\bibitem[Asmussen and Taksar, 1997]{Asmussen1997}
Asmussen, S. and Taksar, M. (1997).
\newblock Controlled diffusion models for optimal dividend pay-out.
\newblock {\em Insurance Math. Econom.}, 20(1):1--15.

\bibitem[Bai and Paulsen, 2010]{Bai-Paulsen-10}
Bai, L. and Paulsen, J. (2010).
\newblock Optimal dividend policies with transaction costs for a class of
  diffusion processes.
\newblock {\em SIAM J. Control Optim.}, 48(8):4987--5008.

\bibitem[Cai et~al., 2009a]{Cai-F-W-09}
Cai, J., Feng, R., and Willmot, G.~E. (2009a).
\newblock Analysis of the compound {P}oisson surplus model with liquid
  reserves, interest and dividends.
\newblock {\em Astin Bull.}, 39(1):225--247.

\bibitem[Cai et~al., 2009b]{Cai-Feng}
Cai, J., Feng, R., and Willmot, G.~E. (2009b).
\newblock On the expectation of total discounted operating costs up to default
  and its applications.
\newblock {\em Adv. in Appl. Probab.}, 41(2):495--522.

\bibitem[Cai et~al., 2009c]{Cai-Feng2}
Cai, J., Feng, R., and Willmot, G.~E. (2009c).
\newblock The compound Poisson surplus model with interest and liquid reserves: analysis of the Gerber-Shiu discounted penalty function.
\newblock {\em Methodology and Computing in Applied Probability} 11: 401-423.

\bibitem[Choulli et~al., 2003]{Choulli}
Choulli, T., Taksar, M., and Zhou, X.~Y. (2003).
\newblock A diffusion model for optimal dividend distribution for a company
  with constraints on risk control.
\newblock {\em SIAM J. Control Optim.}, 41(6):1946--1979.

\bibitem[{de Finetti}, 1957]{DeFi}
{de Finetti}, B. (1957).
\newblock Su un' impostazione alternativa della teoria collettiva del rischio.
\newblock {\em Transactions of the XVth International Congress of Actuaries},
  2:433--443.

\bibitem[Fang and Wu, 2007]{Fang2007}
Fang, Y. and Wu, R. (2007).
\newblock Optimal dividend strategy in the compound {P}oisson model with
  constant interest.
\newblock {\em Stoch. Models}, 23(1):149--166.

\bibitem[Feng et~al., 2012]{pdcp-cdc-2012}
Feng, R., Zhang, S., and Zhu, C. (2012).
\newblock Optimal dividend payment problems in piecewise-deterministic compound
  poisson risk models.
\newblock In {\em Proceeding of the 51st IEEE Conference on Decision and
  Control}, pages 7309--7314, Maui, Hawaii, USA.

\bibitem[Gerber and Shiu, 2006]{Gerber2006}
Gerber, H.~U. and Shiu, E. S.~W. (2006).
\newblock On optimal dividend strategies in the compound {P}oisson model.
\newblock {\em N. Am. Actuar. J.}, 10(2):76--93.

\bibitem[Hartman, 2002]{Hartman}
Hartman, P. (2002).
\newblock {\em Ordinary differential equations}, volume~38 of {\em Classics in
  Applied Mathematics}.
\newblock Society for Industrial and Applied Mathematics (SIAM), Philadelphia,
  PA.
% \newblock Corrected reprint of the second (1982) edition [Birkh{\"a}user,
 % Boston, MA; MR0658490 (83e:34002)], With a foreword by Peter Bates.

\bibitem[Hunting and Paulsen, 2013]{HunPau}
Hunting, M. and Paulsen, J. (2013).
\newblock Optimal dividend policies with transaction costs for a class of jump-diffusion processes.
\newblock {\it Finance and Stochastics}, 17 (1), 73--106.

\bibitem[Jeanblanc-Picqu\'{e} and Shiryaev, 1995]{Jeanblanc1995}
Jeanblanc-Picqu\'{e}, M. and Shiryaev, A.~N. (1995).
\newblock Optimization of the flow of dividends.
\newblock {\em Russian Mathematical Surveys}, 50(2):257--277.

\bibitem[Olver et~al., 2010]{Olver:2010:NHMF}
Olver, F.~W.~J., Lozier, D.~W., Boisvert, R.~F., and Clark, C.~W., editors
  (2010).
\newblock {\em {NIST Handbook of Mathematical Functions}}.
\newblock Cambridge University Press, New York, NY.
\newblock Print companion to \cite{NIST:DLMF}.

\bibitem[Paulsen, 2008]{Paulsen-08}
Paulsen, J. (2008).
\newblock Optimal dividend payments and reinvestments of diffusion processes
  with both fixed and proportional costs.
\newblock {\em SIAM J. Control Optim.}, 47(5):2201--2226.

\bibitem[Schmidli, 2002]{Schmidli-02}
Schmidli, H. (2002).
\newblock On minimizing the ruin probability by investment and reinsurance.
\newblock {\em Ann. Appl. Probab.}, 12(3):890--907.

\bibitem[Schmidli, 2008]{schmidli2008}
Schmidli, H. (2008).
\newblock {\em Stochastic control in insurance}.
\newblock Probability and its Applications (New York). Springer-Verlag London
  Ltd., London.
  
\bibitem[Shreve et~al., 1984]{Shreve-84}
Shreve, S.~E., Lehoczky, J.~P., and Gaver, D.~P. (1984).
\newblock Optimal consumption for general diffusions with absorbing and
  reflecting barriers.
\newblock {\em SIAM J. Control Optim.}, 22(1):55--75.


\bibitem[Song et~al., 2011]{Song-S-Z}
Song, Q.~S., Stockbridge, R.~H., and Zhu, C. (2011).
\newblock On optimal harvesting problems in random environments.
\newblock {\em SIAM J. Control Optim.}, 49(2):859--889.

\bibitem[Yin and Zhu, 2010]{YZ-10}
Yin, G.~G. and Zhu, C. (2010).
\newblock {\em Hybrid Switching Diffusions: Properties and Applications},
  volume~63 of {\em Stochastic Modelling and Applied Probability}.
\newblock Springer, New York.

\bibitem[Zhu, 2011]{Zhu-10}
Zhu, C. (2011).
\newblock Optimal control of risk process in a regime switching environment.
\newblock {\em Automatica}, 47(8):1570--1579.

 
\end{thebibliography}

\def\cprime{$'$}

\vspace{0.1in}

Runhuan Feng

Department of Mathematics

University of Illinois at Urbana-Champaign

Urbana, IL 61801, USA. \qquad  {\tt rfeng@illinois.edu}\\

Hans W. Volkmer

Department of Mathematical Sciences

University of Wisconsin-Milwaukee

Milwaukee, WI 53201, USA. \qquad {\tt volkmer@uwm.edu }\\

Shuaiqi Zhang

School of Science

Hebei University of Technology

Tianjin, 300130, China. \qquad {\tt shuaiqiz@hotmail.com}\\

Chao Zhu

Department of Mathematical Sciences

University of Wisconsin-Milwaukee

Milwaukee, WI 53201, USA. \qquad {\tt zhu@uwm.edu}

\end{document}